\documentclass[11pt]{article}

\usepackage{amsmath,amssymb,amsfonts,amsthm,mathtools,mathrsfs}
\usepackage{authblk}
\usepackage[numbers,sort&compress]{natbib}
\usepackage[a4paper,margin=2.5cm]{geometry}
\usepackage{booktabs,tabularx,array}
\usepackage{microtype}
\usepackage{enumitem}
\usepackage{placeins}
\usepackage[font=small,labelfont=bf]{caption}
\usepackage{titlesec}
\usepackage{needspace}
\usepackage[hyperfootnotes=false,hidelinks]{hyperref}

\hypersetup{
	pdftitle={The maximum adjacency spectral radius of unbalanced Kr-minus-free signed graphs of prescribed size},
	pdfauthor={Chuling Zhang and Xiying Yuan},
	pdfkeywords={signed graph, negative triangle, largest adjacency eigenvalue, prescribed size, spectral extremal graph}
}

\setlist{itemsep=0.25em,topsep=0.5em}

\titleformat{\section}
{\normalfont\normalsize\bfseries}{\thesection.}{0.5em}{}
\titleformat{\subsection}
{\normalfont\normalsize\bfseries}{\thesubsection.}{0.5em}{}
\titlespacing*{\section}{0pt}{2.4ex plus 0.8ex minus 0.2ex}{1.0ex}
\titlespacing*{\subsection}{0pt}{2.0ex plus 0.6ex minus 0.2ex}{0.8ex}

\numberwithin{equation}{section}

\makeatletter
\renewenvironment{proof}[1][\proofname]{%
	\par
	\pushQED{\qed}%
	\normalfont
	\topsep6\p@\@plus6\p@\relax
	\trivlist
	\item[\hskip\labelsep\bfseries #1\@addpunct{.}]%
	\ignorespaces
}{%
	\popQED
	\endtrivlist
	\@endpefalse
}
\makeatother

\theoremstyle{plain}
\newtheorem{theorem}{Theorem}[section]

\newtheorem{lemma}[theorem]{Lemma}

\newtheorem{claim}[theorem]{Claim}
\theoremstyle{definition}
\newtheorem{definition}[theorem]{Definition}
\theoremstyle{remark}

\newcommand{\eps}{\varepsilon}
\newcolumntype{Y}{>{\raggedright\arraybackslash}X}

\title{The maximum adjacency spectral radius of unbalanced
	\(\mathcal K_r^-\)-free signed graphs of prescribed size}

\author[a,b,*]{Chuling Zhang\,$^{\rm a,b}$, \ Xiying Yuan\,}

\affil[a]{\small Department of Mathematics, Shanghai University, Shanghai 200444, PR China}
\affil[b]{\small Newtouch Center for Mathematics of Shanghai University, Shanghai 200444, PR China}

\date{}

\begin{document}
	
	\maketitle
	
	\footnotetext{*\textit{Corresponding author}.}
	\footnotetext{Email addresses: \texttt{zhangchuling727@163.com} (C. Zhang),
		\texttt{xiyingyuan@shu.edu.cn} (X. Yuan).}
	
\begin{abstract}
	Let \(\dot G=(G,\sigma)\) be a signed graph of order \(n\),
	where each edge \(e\in E(G)\) is assigned a weight
	\(\sigma(e)\in\{+1,-1\}\).
	Let \(A(\dot G)\) denote its signed adjacency matrix and
	\(\lambda_1(\dot G)\geq\cdots\geq\lambda_n(\dot G)\)
	be its eigenvalues, counted with multiplicity.
	The largest eigenvalue \(\lambda_1(\dot G)\) is called
	the index of \(\dot G\) and its adjacency spectral radius
	is defined by
	\(\rho(\dot G)=\max\{\lambda_1(\dot G),-\lambda_n(\dot G)\}\).
	Nosal (1970) proved a sharp spectral bound for
	\(C_3\)-free graphs with a prescribed number of edges
	and Nikiforov (2002) extended it to \(K_r\)-free graphs.
	Motivated by these results, we establish signed graph
	analogues by determining the maximum
	index among unbalanced \(\mathcal K_r^-\)-free signed
	graphs of prescribed size \(m\) and characterizing all
	index extremal graphs, where \(\mathcal K_r^-\) denotes
	the family of all unbalanced signed complete graphs
	of order \(r\).
	For the adjacency spectral radius, we determine the
	maximum and the extremal graphs explicitly for \(r=3\)
	and \(m\geq4\).
	For \(r\geq4\), we give explicit determinations
	when \(3\leq m\leq\binom{r}{2}\) or
	\(m\geq \binom{2r-2}{2}+2\), and provide some characterizations in the intermediate
	range \(\binom{r}{2}<m<\binom{2r-2}{2}+2\).
\end{abstract}
	\section{Introduction}
	Throughout, all graphs are finite, simple and undirected.
	A signed graph \(\dot G=(G,\sigma)\) consists of an underlying
	graph \(G\) and a signature \(\sigma:E(G)\to\{-1,+1\}\).
	An edge \(e\) is positive if \(\sigma(e)=+1\) and negative
	if \(\sigma(e)=-1\). A cycle is positive if the product of
	its edge signs is \(+1\), and negative otherwise. A signed
	graph is balanced if it contains no negative cycle, and
	unbalanced otherwise. For \(U\subseteq V(G)\), switching
	with respect to \(U\) reverses the signs of all edges between
	\(U\) and \(V(G)\setminus U\). Two signed graphs on the same
	underlying graph are switching equivalent if one can be
	obtained from the other by a sequence of switchings. They are
	switching isomorphic if one is isomorphic to a signed graph
	obtained from the other by switching.
	Let \(V(G)=\{v_1,\ldots,v_n\}\). The adjacency matrix of
	\(\dot G\) is \(A(\dot G)=(a_{ij}^{\sigma})\), where
	\(a_{ij}^{\sigma}=\sigma(v_iv_j)\) if \(v_iv_j\in E(G)\)
	and \(a_{ij}^{\sigma}=0\) otherwise. 
	Let \(\lambda_1(\dot G)\geq\lambda_2(\dot G)\geq\cdots
	\geq\lambda_n(\dot G)\) be the eigenvalues of \(A(\dot G)\),
	counted with multiplicity. The largest eigenvalue
	\(\lambda_1(\dot G)\) is called the index of \(\dot G\).
	The adjacency spectral radius of \(\dot G\) is defined by
	\(
	\rho(\dot G)
	=\max\{\lambda_1(\dot G),-\lambda_n(\dot G)\}.
	\)
	For an unsigned graph \(G\), we denote its
	adjacency spectral radius by \(\rho(G)\).
	Given a graph \(F\), a graph is called \(F\)-free if it contains
	no subgraph isomorphic to \(F\). For \(r\geq3\), we write
	\(C_r^-\) for a negative \(r\)-cycle and \(\mathcal K_r^-\)
	for the family of unbalanced signed complete graphs of order
	\(r\). A signed graph is \(C_r^-\)-free if it contains no
	negative \(r\)-cycle and \(\mathcal K_r^-\)-free if it contains
	no member of \(\mathcal K_r^-\) as a signed subgraph.
	Throughout this paper, we disregard isolated vertices
	and consider extremal signed graphs up to switching
	isomorphism.

	The spectral Tur\'an problem asks for the maximum adjacency
	spectral radius among \(F\)-free graphs of order \(n\).
	In 2007, Nikiforov~\cite{Nikiforov2007} proved that the
	Tur\'an graph \(T_r(n)\) is the unique extremal graph among
	\(K_{r+1}\)-free graphs of order \(n\). For related results
	on forbidden cycles, see \cite{Nikiforov2009C4,Nikiforov2010}.
	Several results are also known for the Brualdi--Hoffman--Tur\'an
	type problem, which asks for the maximum adjacency spectral radius
	among \(F\)-free graphs with a prescribed number of edges.
	Among all graphs with a given number of edges, Brualdi and
	Hoffman~\cite{BrualdiHoffman1985} conjectured that the
	quasi-complete graph maximizes the adjacency spectral radius.
	Rowlinson~\cite{Rowlinson1988} proved the conjecture and
	established uniqueness up to isolated vertices. For a
	\(C_3\)-free graph \(G\) with \(m\) edges, Nosal~\cite{Nosal1970}
	proved that \(\rho(G)\leq\sqrt{m}\). Nikiforov~\cite{Nikiforov2002}
	extended this result by proving that
	\(\rho(G)\leq\sqrt{2m(1-1/r)}\) for \(K_{r+1}\)-free graphs
	with \(m\) edges. He also proved that \(\rho(G)\leq\sqrt{m}\)
	for \(C_4\)-free graphs with \(m\geq9\), with equality attained
	by a star~\cite{Nikiforov2009C4}. For further results on
	forbidden cycles and complete bipartite subgraphs, see
	\cite{ZhaiLinShu2021,SunLiWei2023} and the survey
	\cite{YuZhangZhang2026}.
	
	Spectral Turán problems and Brualdi--Hoffman--Tur\'an problems have also been studied for signed graphs.
	Kannan and Pragada~\cite{KannanPragada2023} established
	spectral Tur\'an inequalities in terms of the balanced clique
	number. For connected \(C_3^-\)-free unbalanced signed graphs
	of order \(n\), Wang, Hou and Li~\cite{WangHouLi2024}
	determined the maximum adjacency spectral radius and
	characterized the extremal graph. 
	Chen and Yuan~\cite{ChenYuan2024} and
	Wang~\cite{Wang2024K5} studied the corresponding problems
	for \(\mathcal K_4^-\)-free and \(\mathcal K_5^-\)-free
	unbalanced signed graphs, respectively. For \(4\leq r\leq n\),
	Xiong and Hou~\cite{XiongHouCliques2024}
	determined the maximum index among \(\mathcal K_r^-\)-free
	unbalanced signed graphs of order \(n\) and characterized
	the extremal graph up to switching isomorphism.
	For longer negative cycles, Wang and Lin~\cite{WangLin2025}
	determined the maximum index and characterized the extremal
	graph among \(C_4^-\)-free unbalanced signed graphs of order
	\(n\geq5\). Xiong and Hou~\cite{XiongHouCycles2023} obtained
	the corresponding results for \(C_r^-\)-free unbalanced signed
	graphs of order \(n\) when \(5\leq r\leq\lfloor n/3\rfloor+1\).
	Wang, Hou and Huang~\cite{WangHouHuang2025} established a sharp
	adjacency spectral radius bound for the \(C_{2\ell+1}^-\)-free
	case when \(3\leq\ell\leq(n-11)/10\). Cai and Zhou~\cite{CaiZhou2024}
	studied the spectral radius when both \(C_3^-\) and \(C_4^-\)
	are forbidden. For \(n\geq9\), Zha~\cite{Zha2026} characterized
	the unbalanced signed graphs of order \(n\) with maximum index
	among those containing no negative Hamiltonian cycle.
	
	For signed graphs of prescribed size, several works have investigated
	the maximum index under additional constraints. Koledin and
	Stani\'c~\cite{KoledinStanic2017} studied the maximum index of connected
	signed graphs with prescribed order, size and number of negative edges.
	Wang, Peng and Hou~\cite{WangPengHou2026} studied the maximum index of
	unbalanced signed graphs of prescribed size.
	Yang, Li and Chen~\cite{YangLiChen2026} studied the maximum index
	of \(C_4^-\)-free unbalanced signed graphs of prescribed size.
	In this paper, we determine both the maximum index and the maximum
	adjacency spectral radius among \(\mathcal K_r^-\)-free unbalanced signed
	graphs of prescribed size, where \(r\geq3\).

	We first introduce two signed graphs.
	Let \(2\le q\le k\). Consider an all-positive clique \(K_k\)
	with vertex partition
	\(V(K_k)=\{a\}\mathbin{\dot\cup}B\mathbin{\dot\cup}B'\),
	where \(|B|=q-1\) and \(|B'|=k-q.\)
	Here and throughout, \ensuremath{\dot\cup} denotes disjoint union.
	Add two new vertices $u,v$, a negative edge $uv$, and positive edges $ua$
	and $vb$ for every $b\in B$. The resulting
	signed graph is denoted by $\dot\Sigma_{k,q}$. Thus
	\(
	e(\dot\Sigma_{k,q})=\binom{k}{2}+q+1.
	\)
	
	For $d\ge4$, the signed graph $\dot\Omega_d$ is obtained from the
	all-positive graph $K_d-xy$ as follows. Choose two distinct vertices
	$a,b\in V(K_d)\setminus\{x,y\}$. Add a vertex $v$ and a positive edge $vb$,
	and then add a vertex $u$, a positive edge $ua$, and a negative edge
	$uv$. 
	For \(d=3\), begin with the positive path \(xby=K_3-xy\)
	and set \(a=x\). Add a vertex \(v\) and a positive edge
	\(vb\), and then add a vertex \(u\), a positive edge \(ua\),
	and a negative edge \(uv\).
	In either case,
	\(
	e(\dot\Omega_d)=\binom d2+2.
	\)
	
	Clearly, both signed graphs \(\dot\Sigma_{k,q}\) and \(\dot\Omega_d\)
	are \(C_3^-\)-free.
	Moreover, they are unbalanced since
	both contain a \(C_4^-\).
	We can now state our main result.
	
	\begin{theorem}\label{t1.1}
		Let $m\ge4$ and $k$ be the unique integer satisfying
		\(\binom{k}{2}+3\le m\le \binom{k+1}{2}+2.\)
		Let
		\(
		q=m-\binom{k}{2}-1.
		\)
		Then $2\le q\le k+1$.  If $\dot G$ is an unbalanced $C_3^-$-free signed
		graph with $m$ edges, then
		\[
		\lambda_1(\dot G)\le
		\begin{cases}
			\lambda_1(\dot\Sigma_{k,q}),&2\le q\le k,\\[2mm]
			\lambda_1(\dot\Omega_{k+1}),&q=k+1.
		\end{cases}
		\]
		Moreover, equality holds if and only if, after deleting isolated
		vertices, \(\dot G\) is switching isomorphic to the corresponding
		signed graph above.
	\end{theorem}
	 For an unsigned graph \(H\), let \((H,-)\)
	denote the signed graph with every edge of \(H\) negative. 
	The following theorem determines the maximum adjacency spectral
	radius among unbalanced \(C_3^-\)-free signed graphs of prescribed size.
	\begin{theorem}\label{t2}
		Let \(\dot G\) be an unbalanced \(C_3^-\)-free signed graph
		with \(m\geq4\) edges and \(k,q\) be as in
		Theorem~\ref{t1.1}. Then
		\[
		\rho(\dot G)\leq
		\begin{cases}
			2, & m=5,\\[2mm]
			\lambda_1(\dot\Sigma_{k,q}),
			& m\ne5 \text{ and } 2\leq q\leq k,\\[2mm]
			\lambda_1(\dot\Omega_{k+1}),
			& m\ne5 \text{ and } q=k+1.
		\end{cases}
		\]
		Moreover, equality holds if and only if, after deleting
		isolated vertices, \(\dot G\) is switching isomorphic to
		\((C_5,-)\), \(\dot\Sigma_{k,q}\), or \(\dot\Omega_{k+1}\),
		respectively.
	\end{theorem}
	For a nonnegative integer \(p\), let \(R_p\) denote the
	quasi-complete graph with the following construction.
	Let \(R_0\) be the empty graph. If
	\(p=\binom{d}{2}\) for some integer \(d\ge 2\), set
	\(R_p=K_d\). Otherwise, write uniquely
	\(
	p=\binom{d}{2}+t,
	\)
	where \(d\ge 2\) and \(1\le t\le d-1\).
	Then \(R_p\) is obtained from \(K_d\) by adding
	a new vertex adjacent to exactly \(t\) vertices of \(K_d\).
	For the following two theorems, let \(r\geq4\), \(m\geq3\),
	and
	\(
	p=m-2=\binom{k}{2}+s,
	\)
	where \(k\geq2\) and \(0\leq s\leq k-1\) are integers.
	We use two signed graphs with \(m\) edges. To construct \(\dot\Theta_{k,s}\),
	start with the all-positive graph \(R_p\). For \(s=0\), choose distinct
	vertices \(a,b\) of \(K_k\). For \(s>0\), let \(b\) be the vertex outside
	\(K_k\) and choose \(a\in N_{R_p}(b)\). Add a vertex \(u\), a positive edge
	\(ua\), and a negative edge \(ub\). To construct \(\dot B_{k,s+1}\) for
	\(0\leq s\leq k-2\), start with an all-positive \(K_k\) and add a vertex \(u\)
	with \(s+1\) positive neighbors and one negative neighbor in \(K_k\).
	Both graphs contain a negative triangle. In \(\dot\Theta_{k,s}\), an
	unbalanced complete subgraph must contain the negative edge \(ub\), whose
	endpoint \(u\) has degree two. Thus \(\dot\Theta_{k,s}\) is
	\(\mathcal K_4^-\)-free. In \(\dot B_{k,s+1}\), the largest unbalanced
	complete subgraph consists of \(u\) and its \(s+2\) neighbors. Hence
	\(\dot B_{k,s+1}\) is \(\mathcal K_r^-\)-free if and only if \(s\leq r-4\).
	Also, \(\dot B_{k,1}\cong\dot\Theta_{k,0}\).
	\begin{theorem}\label{thm:Kr-index}
		Let \(M=\min\{k-2,r-4\}\). Among the unbalanced \(\mathcal K_r^-\)-free signed
		graphs with \(m\) edges, the graphs attaining the maximum index are as follows.
		\begin{enumerate}[label=\textup{(\roman*)},leftmargin=*,itemsep=3pt]
			\item If \(s=0\) or \(s>M\), the unique extremal graph is
			\(\dot\Theta_{k,s}\).
			
			\item Suppose \(1\leq s\leq M\).
			\begin{enumerate}[label=\textup{(\alph*)},leftmargin=*,itemsep=2pt,topsep=3pt]
				\item If \((k,s)=(3,1)\), the two extremal graphs are
				\(\dot\Theta_{3,1}\) and \(\dot B_{3,2}\), both with index \(\sqrt5\).
				
				\item If \(k\geq3s\) and \((k,s)\ne(3,1)\), the unique extremal graph
				is \(\dot\Theta_{k,s}\).
				
				\item If \(k<3s\), let \(d=3s-k>0\) and \(a=s(k-s-2)\).
				\begin{enumerate}[label=\textup{(\arabic*)},leftmargin=*,itemsep=2pt,topsep=3pt]
					\item If \(\frac{a}{d}\leq k-1\), the unique extremal graph is
					\(\dot B_{k,s+1}\).
					
					\item If \(k-1<\frac{a}{d}<k\), let
					\(
						N=a^3+(2-k)a^2d-(k+s+1)ad^2
						+\bigl[k(s+2)-s^2-s-2\bigr]d^3.
					\)
					Then \(N\ne0\). The unique extremal graph is \(\dot\Theta_{k,s}\)
					for \(N>0\) and \(\dot B_{k,s+1}\) for \(N<0\).
					
					\item If \(\frac{a}{d}\geq k\), the unique extremal graph is
					\(\dot\Theta_{k,s}\).
				\end{enumerate}
			\end{enumerate}
		\end{enumerate}
	\end{theorem}

		For $d\ge3$, let $L_d$ be the graph obtained from $K_d$ by deleting an edge
	$xy$ and attaching a pendant edge at a vertex
	$b\in V(K_d)\setminus\{x,y\}$.  In particular, $L_3=K_{1,3}$.
		Let \(\Lambda_{m,r}\) denote the maximum index in
	Theorem~\ref{thm:Kr-index}.
	Let
	\(
	N_{m,r}
	=\left\lfloor 2m+1-\Lambda_{m,r}^2\right\rfloor
	\)
	and \(\mathcal H_{m,r}\) be the family of connected
	unsigned \(K_r\)-free graphs with \(m\) edges and at most
	\(N_{m,r}\) vertices. Define
	\(
	M_{m,r}
	=\max\left(
	\{\Lambda_{m,r}\}
	\cup
	\{\rho(H):H\in\mathcal H_{m,r}\}
	\right).
	\)
	\begin{theorem}\label{thm:Kr-radius}
		Let \(\dot G\) be an unbalanced \(\mathcal K_r^-\)-free signed graph with
		\(m\) edges and \(T_r=\binom{2r-2}{2}+2\), where \(r\geq4\). Then
		\[
		\rho(\dot G)\leq
		\begin{cases}
			\rho(R_m),      &3\leq m<\binom r2,\\[2pt]
			\rho(L_r),      &m=\binom r2,\\[2pt]
			M_{m,r},        &\binom r2<m<T_r,\\[2pt]
			\Lambda_{m,r},  &m\geq T_r.
		\end{cases}
		\]
	Moreover, equality holds if and only if, after deleting
	isolated vertices, \(\dot G\) is switching isomorphic to
	one of the following signed graphs:
	\begin{enumerate}
		\item[\textup{(i)}]
		\((R_m,-)\), if \(3\le m<\binom{r}{2}\);
		
		\item[\textup{(ii)}]
		\((L_r,-)\), if \(m=\binom{r}{2}\);
		
		\item[\textup{(iii)}]
		if \(\binom{r}{2}<m<T_r\), a maximum index graph
		described in Theorem~\ref{thm:Kr-index} when
		\(\Lambda_{m,r}=M_{m,r}\), or \((H,-)\) for some
		\(H\in\mathcal H_{m,r}\) with \(\rho(H)=M_{m,r}\);
		
		\item[\textup{(iv)}]
		a maximum index graph described in Theorem~\ref{thm:Kr-index},
		if \(m\ge T_r\).
	\end{enumerate}
	\end{theorem}

	\section{Preliminaries}

	We first recall a basic property of switching in signed graphs.
	\begin{lemma}[\cite{SunLiuLan2022}]\label{l1}
		Let \(\dot G\) be a signed graph. Then there exists a signed
		graph \(\dot G'\) switching equivalent to \(\dot G\) such that
		\(A(\dot G')\) has a nonnegative eigenvector corresponding
		to \(\lambda_1(\dot G')\).
	\end{lemma}
	
	The following lemma shows that the index can
	increase strictly even if the Rayleigh quotient remains unchanged.
	It is an
	immediate consequence of the Rayleigh--Ritz theorem; see,
	e.g., \cite[Section 4.2]{HornJohnson2013}.
	\begin{lemma}\label{lem:strict-increase}
		Let \(\dot G\) and \(\dot G'\) be two signed graphs on the
		same vertex set, and \(\mathbf{x}\) be a unit eigenvector
		corresponding to \(\lambda_1(\dot G)\). If
		\(
		\mathbf{x}^{\mathsf T}
		\bigl(A(\dot G')-A(\dot G)\bigr)\mathbf{x}=0
		\)
		and
		\(
		\bigl(A(\dot G')-A(\dot G)\bigr)\mathbf{x}
		\ne \mathbf{0},
		\)
		then
		\(
		\lambda_1(\dot G')>\lambda_1(\dot G).
		\)
	\end{lemma}
	
	\begin{proof}
		By the Rayleigh principle,
		\(
		\lambda_1(\dot G')
		\ge \mathbf{x}^{\mathsf T}A(\dot G')\mathbf{x}
		=\mathbf{x}^{\mathsf T}A(\dot G)\mathbf{x}
		=\lambda_1(\dot G).
		\)
		Suppose that equality holds. Then \(\mathbf{x}\) is also
		an eigenvector corresponding to \(\lambda_1(\dot G')\).
		Consequently,
		\(
		A(\dot G')\mathbf{x}
		=\lambda_1(\dot G')\mathbf{x}
		=\lambda_1(\dot G)\mathbf{x}
		=A(\dot G)\mathbf{x},
		\)
		contradicting
		\(
		\bigl(A(\dot G')-A(\dot G)\bigr)\mathbf{x}\ne\mathbf{0}.
		\)
		Hence \(\lambda_1(\dot G')>\lambda_1(\dot G)\).
	\end{proof}

We next recall the following extremal result for unsigned graphs due to Rowlinson~\cite{Rowlinson1988}.
	
	\begin{lemma}[\cite{Rowlinson1988}]\label{R}
		Let \(G\) be an unsigned graph with \(p\) edges. Then
		\(
		\rho(G)\le \rho(R_p),
		\)
		with equality if and only if
		\(
		G\cong R_p.
		\)
		Moreover,
		\(\rho(R_p)<\rho(R_{p+1})\).
	\end{lemma}
	The strict monotonicity follows from the
	Perron--Frobenius theorem, since \(R_p\) is a proper
	subgraph of the connected graph \(R_{p+1}\) for \(p\ge1\).
	The case \(p=0\) is trivial.
	
	\begin{definition}[\cite{MahadevPeled1995}]\label{def:threshold-graph}
		Let \(G\) be a graph. 
		We call \(G\) a \emph{threshold graph} if it admits a partition
		\(
		V(G)=C\,\dot\cup\,I,
		\)
		where \(C\) is a clique and \(I\) is an independent set
		whose vertices can be labeled \(z_1,\ldots,z_t\) so that
		\(N_G(z_i)\supseteq N_G(z_j)\) for all \(1\leq i<j\leq t\).
	\end{definition}
	\begin{lemma}[{\cite{MahadevPeled1995}}]
		\label{lem:threshold-preorder}
		Let \(G\) be a graph. Define the vicinal preorder on \(V(G)\) by
		\(
		u\preceq v
		\Longleftrightarrow
		N_G(u)\setminus\{v\}\subseteq N_G(v)\setminus\{u\}.
		\)
		Then \(G\) is a threshold graph if
		and only if its vicinal preorder is total, that is, if
		\(u\preceq v\) or \(v\preceq u\) for every \(u,v\in V(G)\).
	\end{lemma}
\begin{theorem}[\cite{WangPengHou2026}]\label{w}
	Let \(d,t\) be positive integers and \(\dot G\)
	be an unbalanced signed graph that maximizes the index in the set
	of unbalanced signed graphs of size \(m\).
	Up to switching isomorphism,
	\begin{enumerate}
		\renewcommand{\labelenumi}{(\arabic{enumi})}
		\item
		if \(m=\binom{d}{2}\), then \(\dot G\) has a unique
		nontrivial component \(\dot B_{d-1,d-2}\);
		
		\item
		if \(m=\binom{d}{2}+1\), then \(\dot G\) has a unique
		nontrivial component \(\dot\Theta_{d-1,d-2}\);
		
		\item
		if \(m=\binom{d}{2}+t+1\), where \(1\le t\le d-2\),
		let
		\(
		\theta=\frac{1}{2}
		\left(3-2d+\sqrt{8d^2-8d+1}\right),
		\)
		we have
		\begin{itemize}
			\item
			if \(t=1\), then \(\dot G\) has a unique
			nontrivial component \(\dot\Theta_{d,0}\), where
			\(\dot\Theta_{d,0}\cong\dot B_{d,1}\);
			
			\item
			if \(2\le t\le\theta-1\), then \(\dot G\) has
			a unique nontrivial component
			\(\dot\Theta_{d,t-1}\);
			
			\item
			if \(\theta\le t\le d-2\), then \(\dot G\) has
			a unique nontrivial component \(\dot B_{d,t}\);
			
			\item
			if \(\theta\) is not an integer and
			\(t=\lfloor\theta\rfloor\), then \(\dot G\)
			has a unique nontrivial component
			\(\dot B_{d,t}\) or \(\dot\Theta_{d,t-1}\).
		\end{itemize}
	\end{enumerate}
\end{theorem}

	\begin{lemma}\label{lem:second-extremal}
		Let \(G\) be an unsigned graph with
		\(\binom{d}{2}\) edges and no isolated vertices, where \(d\geq3\).
		If \(G\not\cong K_d\), then
		\(
		\lambda_1(G)\leq\lambda_1(L_d),
		\)
		with equality if and only if \(G\cong L_d\).
	\end{lemma}
	
	\begin{proof}
		Let
		\(
		\mathcal{G}_d
		=\left\{G: |E(G)|=\binom{d}{2},\ G\not\cong K_d\right\}.
		\)
		We may choose \(G\in\mathcal{G}_d\) such that
		\(
		\lambda_1(G)
		=\max_{G_0\in\mathcal{G}_d}\lambda_1(G_0).
		\)
		Let \(\mathbf{x}\) be a unit eigenvector corresponding to \(\lambda_1(G)\).
		
		We first show that \(G\) is connected. Suppose to the contrary
		that \(G\) is disconnected. Let \(F\) be a component of \(G\) with
		\(\lambda_1(F)=\lambda_1(G)\). Take \(u\in V(F)\) and
		an edge \(ab\) in another component of \(G\). Let
		\(\widetilde{G}=G-ab+ua\) and \(F'\) be the component of \(\widetilde{G}\)
		containing \(F\).
		Since \(ua\) is a bridge of \(\widetilde{G}\) and \(d\geq3\), we have
		\(\widetilde{G}\not\cong K_d\). Moreover, \(|E(\widetilde{G})|=|E(G)|=\binom{d}{2}\).
		Thus \(\widetilde{G}\in\mathcal{G}_d\).
		Since \(F\) and \(F'\) are both connected and \(F\) is a
		proper subgraph of \(F'\), the Perron--Frobenius theorem gives
		\(
		\lambda_1(\widetilde{G})\geq\lambda_1(F')
		>\lambda_1(F)=\lambda_1(G)
		\), a contradiction.

		We consider the following two cases according to whether
		the vicinal preorder of \(G\) is total.

		\medskip
		\noindent
		\textbf{Case 1.} \textit{The vicinal preorder of \(G\) is total.}
		
		By Lemma~\ref{lem:threshold-preorder}, \(G\) is a threshold
		graph. Thus we may write \(V(G)=C\,\dot\cup\,I\), where
		\(C=\{y_1,\ldots,y_c\}\) is a maximum clique and
		\(I=\{z_1,\ldots,z_t\}\) is an independent set. 
		Let
		\(d_i=|N_G(z_i)|\) for \(1\leq i\leq t\). By suitably
		labeling the vertices of \(C\) and \(I\), we may assume that
		\(N_G(z_i)=\{y_1,\ldots,y_{d_i}\}\) and \(c-1\geq d_1\geq\cdots\geq d_t\geq1\), where \(1\leq i\leq t\).
		For convenience, let
		\(\lambda_1=\lambda_1(G)\), \(L=c-1\), and \(D=\sum_{i=1}^{t}d_i
		=\binom{d}{2}-\binom{c}{2}.\)
		
		We may assume that \(x_{y_1}=1\).
		For each \(y_j\in C\), the eigenvalue equations at \(y_1\)
		and \(y_j\) give
		\(
		(\lambda_1+1)(x_{y_1}-x_{y_j})
		=\sum_{\substack{1\leq i\leq t\\ y_j\notin N_G(z_i)}}x_{z_i}
		\geq0.
		\)
		Thus \(x_{y_j}\leq1\) for every \(y_j\in C\). Hence,
		for every \(1\leq i\leq t\),
		\(
		\lambda_1 x_{z_i}=\sum_{j=1}^{d_i}x_{y_j}\leq d_i.
		\)
		The eigenvalue equation
		at \(y_1\) gives
		\(
		\sum_{j=1}^{c}x_{y_j}
		=\lambda_1+1-\sum_{i=1}^{t}x_{z_i}.
		\)
		Summing the eigenvalue equations over all vertices of \(C\),
		we obtain
		\(
		(\lambda_1-L)\sum_{j=1}^{c}x_{y_j}
		=\sum_{i=1}^{t}d_i x_{z_i}.
		\)
		Combining these two equations, we obtain
		\(
		(\lambda_1-L)(\lambda_1+1)
		=\sum_{i=1}^{t}(d_i+\lambda_1-L)x_{z_i}.
		\)
		Since \(G[C]\cong K_c\), we have \(\lambda_1\geq L\).
		Using \(\lambda_1 x_{z_i}\leq d_i\) and
		\(\sum_{i=1}^{t}d_i=D\), we obtain
		\begin{equation}\label{eq:threshold-refined}
			(\lambda_1-L)\lambda_1(\lambda_1+1)
			\leq \sum_{i=1}^{t}d_i^2+(\lambda_1-L)D.
		\end{equation}
		Moreover, since \(d_i\leq L\) for every \(i\), we have
		\(\sum_{i=1}^{t}d_i^2\leq LD\). 
		Substituting this into \eqref{eq:threshold-refined} and
		canceling \(\lambda_1>0\), we obtain
		\begin{equation}\label{eq:threshold-basic}
			(\lambda_1-L)(\lambda_1+1)\leq D.
		\end{equation}
		
		Let \(\alpha=\lambda_1(K_d-e)\), where \(e\) is an edge
		of \(K_d\). By symmetry, a Perron vector of \(K_d-e\)
		has entry \(a\) at each endpoint of \(e\) and entry \(b\)
		at every other vertex. Thus \(\alpha a=(d-2)b\) and
		\(\alpha b=2a+(d-3)b\), which give
		\begin{equation}\label{eq:alpha}
			\alpha^2-(d-3)\alpha-2(d-2)=0.
		\end{equation}
		Since
		\(K_{d-1}\subsetneq K_d-e\subsetneq K_d\)
		and
		\(K_d-e\subsetneq L_d,\)
		the Perron--Frobenius theorem yields
		\(
		d-2
		<\alpha
		<d-1
		\)
		and
		\(
		\alpha<\lambda_1(L_d).
		\)
		By the choice of \(G\), we have
		\(\lambda_1\geq\lambda_1(L_d)>\alpha\).
		Since \(G\) is connected, has \(\binom{d}{2}\) edges,
		and is not isomorphic to \(K_d\), we have \(c\leq d-1\).
		We consider two subcases.
		
		\smallskip
		\noindent
		\textbf{Subcase 1.1.} \textit{\(c\leq d-2\).}
		
		In this subcase, \(L\leq d-3\). For an integer \(\ell\), let
		\(
		\phi(\ell)
		=(\alpha-\ell)(\alpha+1)
		-\binom{d}{2}+\binom{\ell+1}{2}.
		\)
		For \(L\leq\ell<d-3\), since
		\(\ell<d-3<d-2<\alpha\), we have
		\(
		\phi(\ell+1)-\phi(\ell)=\ell-\alpha<0.
		\)
		By \eqref{eq:alpha}, we obtain
		\(
		\phi(L)\geq\phi(d-3)=\alpha-d+2>0.
		\)
		Since \(D=\binom{d}{2}-\binom{L+1}{2}\), it follows that
		\(D<(\alpha-L)(\alpha+1)\). 
		Thus \eqref{eq:threshold-basic} gives
		\(
		(\lambda_1-L)(\lambda_1+1)
		\leq D<(\alpha-L)(\alpha+1).
		\)
		Since \(x\mapsto(x-L)(x+1)\) is strictly increasing for
		\(x\geq L\), we obtain \(\lambda_1<\alpha\), contradicting
		\(\lambda_1>\alpha\).
		
		\smallskip
		\noindent
		\textbf{Subcase 1.2.} \textit{\(c=d-1\).}
		
		Now \(L=d-2\) and \(D=L+1\). By
		\eqref{eq:threshold-refined}, we have
		\(
		\sum_{i=1}^{t}d_i^2
		\geq
		(\lambda_1-L)\bigl(\lambda_1(\lambda_1+1)-(L+1)\bigr).
		\)
		The right-hand side is positive and strictly increasing
		in \(\lambda_1\) for \(\lambda_1>L\).
		Since \(\lambda_1>\alpha>L\), we obtain
		\(
		\sum_{i=1}^{t}d_i^2
		>(\alpha-L)\bigl(\alpha(\alpha+1)-(L+1)\bigr)
		=L^2+L-\alpha
		>L^2-1,
		\)
		where the equality follows from \eqref{eq:alpha} and the
		last inequality follows from \(\alpha<L+1\).
		If \(d_1\leq L-1\), then \(d_i^2\leq(L-1)d_i\) for
		every \(i\), and hence
		\(
		\sum_{i=1}^{t}d_i^2\leq(L-1)D=L^2-1,
		\)
		a contradiction. Therefore \(d_1=L\). Since
		\(\sum_{i=1}^{t}d_i=L+1\) and every \(d_i\) is positive,
		we have \(t=2\) and \(d_2=1\).
		Thus \(C\cup\{z_1\}\) induces \(K_d-y_cz_1\) and \(z_2\)
		is adjacent only to \(y_1\), which is not an endpoint of
		the missing edge. Hence \(G\cong L_d\).

		\medskip
		\noindent
		\textbf{Case 2.} \textit{The vicinal preorder of \(G\) is not total.}
		
		If the vicinal preorder is not total, choose two incomparable
		vertices \(u,v\). Without loss of generality, we may assume
		that \(x_u\geq x_v\). Then
		\(
		\mathcal{S}=N_G(v)\setminus\bigl(N_G(u)\cup\{u\}\bigr)
		\neq\varnothing.
		\)
		Let \(G'\) be obtained from \(G\) by replacing \(vz\) with
		\(uz\) for every \(z\in\mathcal{S}\). Then
		\(
		\mathbf{x}^{\mathsf T}
		\bigl(A(G')-A(G)\bigr)\mathbf{x}
		=2(x_u-x_v)\sum_{z\in\mathcal{S}}x_z\geq0.
		\)
		If \(x_u>x_v\), the Rayleigh principle gives
		\(\lambda_1(G')>\lambda_1(G)\). If \(x_u=x_v\), then
		\(
		\bigl((A(G')-A(G))\mathbf{x}\bigr)_u
		=\sum_{z\in\mathcal{S}}x_z>0,
		\)
		and by Lemma~\ref{lem:strict-increase}, we also have \(\lambda_1(G')>\lambda_1(G)\). Since the transformation preserves the number
		of edges, the maximality of \(G\) implies that \(K_d\) is
		the only nontrivial component of \(G'\). 
		Since \(uz\in E(G')\) and \(vz\notin E(G')\) for any
		\(z\in\mathcal{S}\), the vertex \(v\) lies outside the
		\(K_d\) component and is therefore isolated in \(G'\).
		Thus \(uv\notin E(G)\) and
		\(N_G(u)\cap N_G(v)=\varnothing\).
		Let \(W=V(G)\setminus\{u,v\}\). Then we have
		\(G[W]=G'[W]\cong K_{d-1}\). 
		Since \(u\) is adjacent
		to every vertex of \(W\) in \(G'\), we have
		\(
		W=N_G(u)\,\dot\cup\,N_G(v).
		\)
		
		Choose \(z\in N_G(v)\). Let \(M_d\) be obtained from
		\(G\) by replacing \(vw\) with \(uw\) for every
		\(w\in N_G(v)\setminus\{z\}\). Thus \(M_d\) consists of
		\(K_d-uz\) on \(W\cup\{u\}\) and the pendant edge \(vz\).
		By the Rayleigh principle,
		\(
		\lambda_1(M_d)-\lambda_1(G)
		\geq \mathbf{x}^{\mathsf T}
		\bigl(A(M_d)-A(G)\bigr)\mathbf{x}
		=2(x_u-x_v)
		\sum_{w\in N_G(v)\setminus\{z\}}x_w
		\geq0.
		\)
		Thus \(\lambda_1(M_d)\geq\lambda_1(G)\).
		Let \(\mathbf{y}\) be a
		positive unit eigenvector of \(A(M_d)\) corresponding to
		\(\lambda_1(M_d)\). Choose \(b\in W\setminus\{z\}\).
		Subtracting
		the eigenvalue equations at \(u\) and \(v\), and those
		at \(b\) and \(z\), we obtain
		\(\lambda_1(M_d)(y_u-y_v)=\sum_{w\in W\setminus\{z\}}y_w-y_z\)
		and
		\((\lambda_1(M_d)+1)(y_b-y_z)=y_u-y_v.\)
		Combining these two equations gives
		\(
		\bigl(\lambda_1(M_d)(\lambda_1(M_d)+1)-1\bigr)(y_b-y_z)
		=\sum_{w\in W\setminus\{z,b\}}y_w\geq0.
		\)
		Since \(M_d\) contains an edge, we have \(\lambda_1(M_d)\geq1\),
		and hence \(y_b\geq y_z\).
		
		Let \(M_d'\) be obtained from \(M_d\) by replacing \(vz\)
		with \(vb\). Since \(b\) is not an endpoint of the missing
		edge \(uz\), we have \(M_d'\cong L_d\). Moreover,
		\(
		\lambda_1(M_d')-\lambda_1(M_d)
		\geq \mathbf{y}^{\mathsf T}
		\bigl(A(M_d')-A(M_d)\bigr)\mathbf{y}
		=2y_v(y_b-y_z)\geq0.
		\)
		If \(y_b>y_z\), the Rayleigh principle gives
		\(\lambda_1(M_d')>\lambda_1(M_d)\). If \(y_b=y_z\), then
		\(
		\bigl((A(M_d')-A(M_d))\mathbf{y}\bigr)_b=y_v>0,
		\)
		and Lemma~\ref{lem:strict-increase} also gives \(\lambda_1(M_d')>\lambda_1(M_d)\). 
		Consequently,
		\(
		\lambda_1(L_d)>\lambda_1(M_d)\geq\lambda_1(G),
		\)
		contrary to the choice of \(G\). Thus Case~2 cannot occur,
		and Case~1 gives \(G\cong L_d\).
	\end{proof}

	\begin{lemma}\label{lem:one-component}
		Let \(r\ge3\) and \(\dot G=(G,\sigma)\) be an unbalanced
		\(\mathcal K_r^-\)-free signed graph with \(m\) edges. If \(\dot G\) has
		maximum index among all unbalanced \(\mathcal K_r^-\)-free signed graphs
		with \(m\) edges, then \(\dot G\) is connected.
	\end{lemma}
	
	\begin{proof}
		Suppose to the contrary that \(\dot G\) is disconnected.
		Then \(\dot G\) has at least two connected components.
		Let \(\dot G_1\) be a connected component of \(\dot G\)
		with \(\lambda_1(\dot G)=\lambda_1(\dot G_1)\).
		Let \(\mathbf{x}\) be a unit eigenvector of
		\(A(\dot G_1)\) corresponding to \(\lambda_1(\dot G_1)\).
		Define \(\widetilde{\mathbf{x}}\in\mathbb{R}^{V(G)}\) by
		\[
		\widetilde{x}_v=
		\begin{cases}
			x_v, & v\in V(G_1),\\
			0,   & v\in V(G)\setminus V(G_1).
		\end{cases}
		\]
		Then \(\widetilde{\mathbf{x}}\) is a unit eigenvector of
		\(A(\dot G)\) corresponding to \(\lambda_1(\dot G)\).
		Choose \(z\in V(G_1)\) with \(\widetilde{x}_z\ne0\).
		Let \(\dot G_2\) be another connected component of \(\dot G\)
		and choose \(w\in V(G_2)\). Since \(\dot G\) has no isolated
		vertices, \(N_G(w)\ne\varnothing\). 
		Let \(\dot G'\) be obtained from \(\dot G\) by replacing the edge
		\(wv\) with an edge \(zv\) of sign \(\sigma(wv)\) for each
		\(v\in N_G(w)\).
		Clearly, \(\dot G'\) still has \(m\) edges.

		We first prove that \(\dot G'\) is unbalanced. Let \(C\) be a
		negative cycle of \(\dot G\). If \(C\) does not contain \(w\),
		then it remains a negative cycle in \(\dot G'\). Otherwise,
		write \(C=wv_1\cdots v_tw\). Since \(z\notin V(G_2)\), replacing
		\(w\) with \(z\) gives a cycle \(zv_1\cdots v_tz\) in \(\dot G'\).
		The edges \(zv_1\) and \(v_tz\) have the same signs as \(wv_1\)
		and \(v_tw\), respectively, and all other edges of the cycle
		are unchanged. Hence this cycle is negative.
		
		We next claim that \(\dot G'\) is \(\mathcal K_r^-\)-free.
		Suppose to the contrary that \(\dot G'\) contains a member of \(\mathcal K_r^-\).
		Since \(\dot G\) is \(\mathcal K_r^-\)-free, this \(\mathcal K_r^-\) must contain
		a new edge \(zv\). In \(\dot G'\), the only neighbor of \(v\) outside
		\(V(G_2)\) is \(z\). Thus all vertices of
		this \(\mathcal K_r^-\) other than \(z\) belong to
		\(V(G_2)\setminus\{w\}\). 
		Replacing \(z\) with \(w\) shows that \(\dot G\) contains
		an unbalanced \(K_r\), a contradiction.
		Thus \(\dot G'\) is \(\mathcal K_r^-\)-free.
		Since \(\widetilde{x}_w=\widetilde{x}_v=0\) for every
		\(v\in N_G(w)\), we have
		\(
		\widetilde{\mathbf{x}}^{\mathsf T}
		\bigl(A(\dot G')-A(\dot G)\bigr)\widetilde{\mathbf{x}}
		=2\sum_{v\in N_G(w)}
		\sigma(wv)(\widetilde{x}_z-\widetilde{x}_w)\widetilde{x}_v
		=0.
		\)
		On the other hand, for any \(v\in N_G(w)\),
		\(
		\bigl((A(\dot G')-A(\dot G))\widetilde{\mathbf{x}}\bigr)_v
		=\sigma(wv)(\widetilde{x}_z-\widetilde{x}_w)
		=\sigma(wv)\widetilde{x}_z\ne0.
		\)
		By Lemma~\ref{lem:strict-increase},
		\(\lambda_1(\dot G')>\lambda_1(\dot G)\), contradicting the
		extremality of \(\dot G\). Hence \(\dot G\) is connected.
	\end{proof}
	
	\begin{lemma}\label{lem:one-negative-edge}
		Let \(r\ge3\) and \(\dot G\) be an unbalanced \(\mathcal K_r^-\)-free signed graph
		with \(m\) edges and maximum index among all such graphs.
		Then \(\dot G\) is switching equivalent to a signed graph
		\(\dot G^\ast\) with exactly one negative edge and
		\(A(\dot G^\ast)\) has a nonnegative unit eigenvector
		corresponding to \(\lambda_1(\dot G^\ast)\).
	\end{lemma}
	
	\begin{proof}
		By Lemma~\ref{lem:one-component}, \(\dot G\) is connected.
		By Lemma~\ref{l1}, there exists a signed graph \(\dot G^\ast\)
		switching equivalent to \(\dot G\) such that \(A(\dot G^\ast)\)
		has a nonnegative unit eigenvector \(\mathbf{x}\) corresponding
		to \(\lambda_1(\dot G^\ast)\). Since switching preserves the
		underlying graph, the index, and the sign of every cycle,
		\(\dot G^\ast\) is also extremal.
		Let \(H\) be the spanning subgraph of \(G\) consisting of
		the positive edges of \(\dot G^\ast\). Denote \(r^-=|E^-(\dot G^\ast)|\),
		where \(E^-(\dot G^\ast)\) is the set of negative edges of
		\(\dot G^\ast\). Then \(e(H)=m-r^-\). Since
		\(\mathbf{x}\geq\mathbf{0}\), the Rayleigh principle gives
		\begin{equation}\label{eq:positive-subgraph-bound}
			\begin{aligned}
				\lambda_1(\dot G^\ast)
				=\mathbf{x}^{\mathsf T}A(H)\mathbf{x}
				-2\sum_{uv\in E^-(\dot G^\ast)}x_u x_v
				\leq \mathbf{x}^{\mathsf T}A(H)\mathbf{x}
				\leq \lambda_1(H).
			\end{aligned}
		\end{equation}
		
		Since \(\dot G^\ast\) is unbalanced, \(r^-\geq1\).
		We claim that \(\lambda_1(\dot G^\ast)<\lambda_1(H)\).
		If \(\lambda_1(\dot G^\ast)=\lambda_1(H)\), then equality
		throughout \eqref{eq:positive-subgraph-bound} gives
		\(A(H)\mathbf{x}=\lambda_1(\dot G^\ast)\mathbf{x}\)
		and
		\(x_u x_v=0\) for every \(uv\in E^-(\dot G^\ast).\)
		Since \(A(G)=2A(H)-A(\dot G^\ast)\), we have
		\(A(G)\mathbf{x}=\lambda_1(\dot G^\ast)\mathbf{x}\), and hence
		\(\sum_{u\in N_G(v)}x_u=\lambda_1(\dot G^\ast)x_v\) for every \(v\in V(G).\)
		If \(x_v=0\) for some \(v\in V(G)\), then the sum is zero.
		Since \(\mathbf{x}\geq\mathbf{0}\), we have \(x_u=0\) for every
		\(u\in N_G(v)\). Since \(G\) is connected, applying the same argument
		along paths from \(v\) gives \(\mathbf{x}=\mathbf{0}\), a contradiction
		to \(\|\mathbf{x}\|_2=1\). Hence \(x_v>0\) for every \(v\in V(G)\),
		contradicting
		\(x_u x_v=0\) on each negative edge. Thus
		\(\lambda_1(\dot G^\ast)<\lambda_1(H)\).
		
		We claim that \(r^-=1\). Suppose to the contrary that \(r^-\geq2\). If \(r^-=2\),
		a negative cycle of \(\dot G^\ast\) contains exactly one
		negative edge, say \(ab\). Deleting \(ab\) from this cycle
		leaves an \(a\)--\(b\) path in \(H\). Let \(\dot G'\) be
		obtained from \(H\) by adding a new vertex \(w\), a positive
		edge \(aw\), and a negative edge \(wb\). The path in \(H\)
		together with \(awb\) forms a negative cycle. 
		Since \(ab\notin E(H)\), \(\dot G'\) is \(\mathcal K_r^-\)-free.
		Hence, \(\dot G'\) is unbalanced and
		\(\mathcal K_r^-\)-free, has \(m\) edges, and contains \(H\) as an
		all-positive induced subgraph. By eigenvalue interlacing,
		\(
		\lambda_1(\dot G')\geq\lambda_1(H)>\lambda_1(\dot G^\ast),
		\)
		contradicting the extremality of \(\dot G^\ast\).
		If \(r^-\geq3\), choose \(ab\in E(H)\), which exists since
		\(\lambda_1(H)>\lambda_1(\dot G^\ast)\geq0\). Let \(\dot G''\)
		be obtained from \(H\) by adding a path \(auvb\) with two
		new vertices \(u,v\), where \(uv\) is negative and \(au,vb\)
		are positive, and attaching \(r^--3\) positive pendant
		edges at \(a\). Then \(auvba\) is a negative cycle.
		Since \(uv\) is the only negative edge and
		\(N_{G''}(u)\cap N_{G''}(v)=\varnothing\), \(\dot G''\) is \(\mathcal K_r^-\)-free.
		Since \(\dot G''\) has \(m\) edges and contains \(H\) as an
		all-positive induced subgraph, the same interlacing argument gives
		\(
		\lambda_1(\dot G'')\geq\lambda_1(H)>\lambda_1(\dot G^\ast),
		\)
		contradicting the extremality of \(\dot G^\ast\). Hence \(r^-=1\).
	\end{proof}

	\section{Spectral extremal results for unbalanced
		\(\mathcal K_3^-\)-free signed graphs}
	
	\subsection{A complete characterization of the index extremal graphs}
	
	We next record the spectral facts needed for the two candidates
	in Theorem~\ref{t1.1}.
	
	\begin{lemma}\label{lem:candidate-lower}
		Let \(p=m-2\) and \(k,q\) be as in
		Theorem~\ref{t1.1}.
		\begin{enumerate}[label=\textup{(\roman*)}]
			\item If \(2\leq q\leq k\), then
			\(
			\lambda_1(\dot\Sigma_{k,q})\geq\lambda_1(R_p).
			\)
			\item If \(q=k+1\) and \(d=k+1\geq4\), then
			\(\tau=\lambda_1(\dot\Omega_d)\) is the unique positive root of
			\[
			t^3-(d-3)t^2-(2d-4)t-2=0
			\]
			and \(d-2<\tau<d-1\).
		\end{enumerate}
	\end{lemma}
	
	\begin{proof}
		We first prove \textup{(i)}. Since
		\(p=\binom{k}{2}+q-1\), deleting \(u\) from
		\(\dot\Sigma_{k,q}\) leaves the all-positive graph \(R_p\).
		Let \(\mathbf{x}\) be a positive unit eigenvector of
		\(A(R_p)\) corresponding to \(\lambda_1(R_p)\) and 
		\(\widehat{\mathbf{x}}\) be obtained from \(\mathbf{x}\)
		by assigning zero to the coordinate at \(u\).
		By the Rayleigh principle, we have
		\(
		\lambda_1(\dot\Sigma_{k,q})
		\geq
		\widehat{\mathbf{x}}^{\mathsf T}
		A(\dot\Sigma_{k,q})\widehat{\mathbf{x}}
		=\mathbf{x}^{\mathsf T}A(R_p)\mathbf{x}
		=\lambda_1(R_p).
		\)
		
		We next prove \textup{(ii)}. Let \(\tau=\lambda_1(\dot\Omega_d)\).
		Since \(V(K_d)\setminus\{x\}\) induces an all-positive \(K_{d-1}\)
		in \(\dot\Omega_d\), interlacing gives \(\tau\geq d-2\geq2\).
		Let \(\mathbf{z}\) be an eigenvector corresponding to \(\tau\)
		and \(S=\sum_{i\in V(K_d)}z_i\). 
		The eigenvalue equations at \(u,v,a,b\) give
		\((\tau-1)(z_u-z_v)=z_a-z_b\)
		and \((\tau+1)(z_a-z_b)=z_u-z_v.\)
		Combining these identities gives
		\((\tau^2-2)(z_u-z_v)=0\).
		Since \(\tau\geq2\), we have \(z_u=z_v\), and hence \(z_a=z_b\).
		The eigenvalue equations at \(u\) and \(a\)
		give \(z_a=(\tau+1)z_u\) and
		\((\tau+1)z_a=S+z_u\), respectively. Hence
		\(	z_u=z_v=\frac{S}{\tau(\tau+2)}\)
		and
		\(z_a=z_b=\frac{(\tau+1)S}{\tau(\tau+2)}.\)
		The equations at \(x\) and \(y\) give
		\(\tau z_x=\tau z_y=S-z_x-z_y\), so \(z_x=z_y=\frac{S}{\tau+2}\).
		For every \(i\in V(K_d)\setminus\{a,b,x,y\}\), we have
		\(\tau z_i=S-z_i\), and hence \(z_i=\frac{S}{\tau+1}\).
		If \(S=0\), these formulas give \(\mathbf{z}=\mathbf{0}\),
		a contradiction. Thus \(S\ne0\). Summing the coordinates over
		\(V(K_d)\) and dividing by \(S\), we obtain
		\(
		1=\frac{2(\tau+1)}{\tau(\tau+2)}
		+\frac{2}{\tau+2}+\frac{d-4}{\tau+1}.
		\)
		This equation is equivalent to \(f_d(\tau)=0\), where
		\(
		f_d(t)=t^3-(d-3)t^2-(2d-4)t-2.
		\)
		By Descartes' rule of signs, \(f_d\) has at most one positive
		root. Since \(\tau>0\), it is the unique positive root. Moreover,
		\(
		f_d(d-2)=-(d-2)^2-2<0<2(d-2)=f_d(d-1).
		\)
		Therefore \(d-2<\tau<d-1\).
	\end{proof}
	
	Let $H$ be a connected unsigned graph and $a,b$ be
	nonadjacent vertices of $H$. Denote by $T(H;a,b)$ the signed
	graph obtained from $H$ by adding a new vertex $w$, a positive
	edge $wa$, and a negative edge $wb$. A positive $a$--$b$ path in $H$ together with $wa$ and $wb$
	forms a negative cycle.
	Thus $T(H;a,b)$ is unbalanced. Since $ab\notin E(H)$, the signed graph is also $C_3^-$-free.

	\begin{lemma}\label{lem:rooted-threshold}
		Suppose that $\Gamma=T(H;a,b)$ maximizes $\lambda_1$
		among all unbalanced $C_3^-$-free signed graphs with
		$|E(\Gamma)|$ edges. Let $\mathbf{x}\geq0$ be a unit
		eigenvector of $A(\Gamma)$ corresponding to
		$\lambda_1(\Gamma)$. Suppose that
		\(x_w\leq x_b,\)
		\(x_b>0,\)
		and \(\lambda_1(\Gamma)>1.\)
		Then $H$ is a threshold graph. More precisely, it admits
		a partition
		\(V(H)=C\mathbin{\dot\cup}I\), \(a\in C\), \(b\in I\),
		where $C$ is a clique and $I$ is an independent set.
		The vertices can be labelled as
		\(C=\{y_1,\ldots,y_c\}\) and
		\(I=\{z_1,\ldots,z_t=b\},\)
		so that
		\(N_H(z_i)=\{y_1,\ldots,y_{d_i}\}\)
		for all \(1\leq i\leq t\)
		and \(c-1\geq d_1\geq\cdots\geq d_t=h\geq1.\)
		Moreover, $y_1\ne a$ and
		\(
		x_{y_1}=\max_{y\in C}x_y.
		\)
	\end{lemma}
	
	\begin{proof}
		We call a signed graph \emph{feasible} if it has the same
		number of edges as $\Gamma$, is unbalanced, and is $C_3^-$-free.
		The eigenvalue equation at $w$ gives \(\lambda_1(\Gamma) x_w=x_a-x_b.\)
		Thus $x_a\geq x_b>0$. We claim that $x_v>0$ for every
		$v\in V(H)$. Suppose that $x_v=0$ for some
		$v\in V(H)\setminus\{b\}$. Since all edges incident with $v$
		in $\Gamma$ are positive and $\mathbf{x}\geq0$, the eigenvalue
		equation at $v$ gives $x_u=0$ for every $u\in N_H(v)$.
		As $H$ is connected, applying the same argument along a
		$v$--$b$ path in $H$ gives $x_b=0$, a contradiction.

		We first show that $d_H(a,b)=2$. Since $ab\notin E(H)$,
		we have $d_H(a,b)\geq2$. Suppose that $d_H(a,b)=\ell\geq3$
		and let $a=v_0,v_1,\ldots,v_{\ell-1},v_\ell=b$ be a shortest
		$a$--$b$ path in $H$. Then $bv_1,av_{\ell-1}\notin E(H)$.
		If $x_{v_1}\geq x_{v_{\ell-1}}$, let $\Gamma_1$ be obtained
		from $\Gamma$ by replacing the positive edge $bv_{\ell-1}$
		with $bv_1$. Then $\Gamma_1$ is feasible. By the Rayleigh principle, we have
		\(
		\lambda_1(\Gamma_1)-\lambda_1(\Gamma)
		\geq 2x_b(x_{v_1}-x_{v_{\ell-1}})\geq0.
		\)
		If $x_{v_1}>x_{v_{\ell-1}}$, we have
		$\lambda_1(\Gamma_1)>\lambda_1(\Gamma)$, contradicting the extremality
		of $\Gamma$. 
		If $x_{v_1}=x_{v_{\ell-1}}$, then
		\(
		\mathbf{x}^{\mathsf T}
		\bigl(A(\Gamma_1)-A(\Gamma)\bigr)\mathbf{x}=0,
		\)
		but
		\(
		\bigl((A(\Gamma_1)-A(\Gamma))\mathbf{x}\bigr)_{v_1}
		=x_b>0.
		\)
		Thus Lemma~\ref{lem:strict-increase} also gives
		$\lambda_1(\Gamma_1)>\lambda_1(\Gamma)$, a contradiction.
		If $x_{v_1}<x_{v_{\ell-1}}$, replace the positive edge $av_1$
		with $av_{\ell-1}$. A similar argument gives a contradiction.
		Hence $d_H(a,b)=2$, and therefore
		$N_H(a)\cap N_H(b)\ne\varnothing$.

		We next show that $N_H(b)\subseteq N_H(a)$. Suppose that
		$y\in N_H(b)\setminus N_H(a)$ and replace the positive edge
		$by$ with $ay$. Clearly, the resulting graph is feasible. Since $x_a\geq x_b$, the same argument as
		above gives a contradiction. Hence $N_H(b)\subseteq N_H(a)$.
		
		We claim that for every $z\in V(H)\setminus\{a,b\}$,
		$x_a\geq x_z$ if $zb\notin E(H)$ and $x_b\leq x_z$ if
		$az\notin E(H)$. If $zb\notin E(H)$, replacing the positive
		edge $wa$ with the positive edge $wz$ gives the feasible graph
		$T(H;z,b)$. If $az\notin E(H)$, replacing the negative edge
		$wb$ with the negative edge $wz$ gives the feasible graph
		$T(H;a,z)$. The same argument as above proves both inequalities.

		We next show that $x_y\geq x_a$ for every $y\in N_H(b)$.
		Suppose that $x_y<x_a$ for some $y\in N_H(b)$. Let $\Gamma_2$
		be obtained from $\Gamma$ by replacing the positive edges
		$wa$ and $by$ with $wy$ and $ab$, respectively. 
		Since $N_H(b)\subseteq N_H(a)$, $\Gamma_2$ contains
		the negative cycle $wyabw$.
		Moreover, $y\nsim b$ in $\Gamma_2$, so $\Gamma_2$ is feasible.
		By the Rayleigh principle, we have
		\(
		\lambda_1(\Gamma_2)-\lambda_1(\Gamma)
		\geq 2(x_y-x_a)(x_w-x_b)\geq0.
		\)
		If $x_w<x_b$, the right-hand side is positive, contradicting
		the extremality of $\Gamma$. If $x_w=x_b$, then
		\(
		\mathbf{x}^{\mathsf T}
		\bigl(A(\Gamma_2)-A(\Gamma)\bigr)\mathbf{x}=0,
		\)
		but
		\(
		\bigl((A(\Gamma_2)-A(\Gamma))\mathbf{x}\bigr)_w
		=x_y-x_a\ne0.
		\)
		Thus Lemma~\ref{lem:strict-increase} also gives a contradiction.
		Hence $x_y\geq x_a$ for every $y\in N_H(b)$.
		
		We next compare the vertices in the vicinal preorder of $H$.
		For distinct vertices $u,v\in V(H)$, let
		\(
		S=N_H(v)\setminus\bigl(N_H(u)\cup\{u\}\bigr).
		\)
		The Kelmans transformation from $v$ to $u$ replaces the
		positive edge $vr$ with $ur$ for every $r\in S$. If
		$x_u\geq x_v$ and the resulting graph $\Gamma'$ is feasible,
		then $v\preceq u$. Indeed, if $S\ne\varnothing$, the change
		in the Rayleigh quotient is
		\(
		2(x_u-x_v)\sum_{r\in S}x_r\geq0
		\)
		and
		\(
		\bigl((A(\Gamma')-A(\Gamma))\mathbf{x}\bigr)_u
		=\sum_{r\in S}x_r>0.
		\)
		The same argument as above gives a contradiction. Thus
		$S=\varnothing$, which means that $v\preceq u$.
		
		For $y\in N_H(b)$, the transformation from $a$ to $y$
		preserves the negative cycle $waybw$ and $a\nsim b$ still holds.
		Since $x_y\geq x_a$, we obtain $a\preceq y$. 
		For $z\in V(H)\setminus\bigl(N_H(b)\cup\{a,b\}\bigr)$,
		the Kelmans transformation from $z$ to $a$ clearly yields
		a feasible graph. Since $x_a\geq x_z$, we obtain $z\preceq a$. 
		Also, $N_H(b)\subseteq N_H(a)$
		and $ab\notin E(H)$ give $b\preceq a$.
		
		We claim that $b\preceq z$ for every
		$z\in V(H)\setminus\{a,b\}$. If $z\in N_H(b)$, this follows
		from $b\preceq a\preceq z$. If $z\notin N_H(b)$ and
		$az\in E(H)$, then $a\preceq y$ implies $z\in N_H(y)$ for
		every $y\in N_H(b)$. Thus $N_H(b)\subseteq N_H(z)$, and
		hence $b\preceq z$.
		It remains to consider the case $az,bz\notin E(H)$.
		Since $az\notin E(H)$, we have $x_z\geq x_b$. Suppose that
		$S'=N_H(b)\setminus N_H(z)$ is nonempty. 
		Let $\Gamma_3$ be
		obtained from $\Gamma$ by replacing the negative edge $wb$
		with the negative edge $wz$ and the positive edge $by$ with the positive edge $zy$ for every
		$y\in S'$. 
		Then $\Gamma_3$ is feasible, since it contains the
		negative cycle $wayzw$ for any $y\in S'$ and $a\nsim z$ still holds in $\Gamma_3$.
		For every $y\in S'\subseteq N_H(b)$, we have
		\(
		x_y\geq x_a=x_b+\lambda_1(\Gamma)x_w>x_w.
		\)
		By the Rayleigh principle, we obtain
		\(
		\lambda_1(\Gamma_3)-\lambda_1(\Gamma)
		\geq 2(x_z-x_b)\left(\sum_{y\in S'}x_y-x_w\right)\geq0.
		\)
		If $x_z>x_b$, the right-hand side is positive, a contradiction.
		If $x_z=x_b$, the Rayleigh quotient is unchanged, but
		\(
		\bigl((A(\Gamma_3)-A(\Gamma))\mathbf{x}\bigr)_b
		=x_w-\sum_{y\in S'}x_y<0.
		\)
		Thus Lemma~\ref{lem:strict-increase} again gives a contradiction.
		Hence $S'=\varnothing$, so $N_H(b)\subseteq N_H(z)$ and
		$b\preceq z$. 
		
		We claim that $H$ is a threshold graph. By
		Lemma~\ref{lem:threshold-preorder} and the preceding comparisons,
		it suffices to compare distinct $u,v\in V(H)\setminus\{a,b\}$.
		Without loss of generality, assume that $x_u\geq x_v$.
		If $u$ and $v$ were incomparable, there would be a vertex
		$z\in N_H(v)\setminus\bigl(N_H(u)\cup\{u\}\bigr)$.
		The Kelmans transformation from $v$ to $u$ then transfers
		at least one edge and yields a feasible graph, since $a\nsim b$
		still holds and $a$ and $b$ still have a common neighbor.
		The same argument as above gives a contradiction.
		Hence $H$ is a threshold graph.

		Since $H$ is a threshold graph, it admits a split partition
		$V(H)=C_0\mathbin{\dot\cup}I_0$, where $C_0$ is a maximum
		clique and $I_0$ is an independent set. We claim that such a
		partition $V(H)=C\mathbin{\dot\cup}I$ can be chosen with
		$a\in C$ and $b\in I$. If $a\in C_0$, take $C=C_0$ and
		$I=I_0$. Then $b\in I$ since $a\nsim b$. Now suppose that
		$a\notin C_0$, so $a\in I_0$.
		Since $C_0$ is a maximum clique, there exists
		$c_0\in C_0\setminus N_H(a)$. Since $N_H(b)\subseteq N_H(a)$,
		we have $c_0\notin N_H(b)$. The preceding comparisons give
		$c_0\preceq a$, including the case $c_0=b$. Since $c_0\nsim a$,
		we have $N_H(c_0)\subseteq N_H(a)$. Hence
		\(
		C=(C_0\setminus\{c_0\})\cup\{a\}
		\)
		is a maximum clique. 
		Moreover, $a$ has no neighbor in
		$I_0\setminus\{a\}$, so neither does $c_0$. Therefore
		\(
		I=(I_0\setminus\{a\})\cup\{c_0\}
		\)
		is an independent set. 
		Moreover,
		$b\preceq z$ for every $z\in I\setminus\{b\}$. We may therefore
		label $I=\{z_1,\ldots,z_t\}$ so that $z_t=b$ and
		\(
		N_H(z_1)\supseteq\cdots\supseteq N_H(z_t)=N_H(b).
		\)
		Since $b\in I$ and $I$ is independent, we have $N_H(b)\subseteq C$.
		Choose $y_1\in N_H(b)$. Since $y_1$ belongs to every $N_H(z_i)$,
		we can label the remaining vertices of $C$ so that
		\(N_H(z_i)=\{y_1,\ldots,y_{d_i}\}\)
		and
		\(c-1\geq d_1\geq\cdots\geq d_t=h\geq1,\)
		where \(1\leq i\leq t,\) $c=|C|$, $d_i=|N_H(z_i)|$, and $h=|N_H(b)|$.
		Since $C$ is a maximum clique, we have $d_1\leq c-1$.
		Since $H$ is connected,
		$h=|N_H(b)|\geq1$.
		
		Finally, $y_1\ne a$ since $a\nsim b$. For every
		$y\in C\setminus\{a,y_1\}$, the eigenvalue equations at
		$y_1$ and $y$ give
		\(
		\bigl(\lambda_1(\Gamma)+1\bigr)(x_{y_1}-x_y)
		=\sum_{\substack{1\leq i\leq t\\ y\notin N_H(z_i)}}x_{z_i}
		\geq0.
		\)
		Thus \(x_{y_1}\geq x_y\) for every
		$y\in C\setminus\{a,y_1\}$.
		Since $y_1\in N_H(b)$, we also have $x_{y_1}\geq x_a$.
		Hence
		$x_{y_1}=\max_{y\in C}x_y$, as required.
	\end{proof}

	\begin{lemma}\label{lem:split-ear}
		Retain the notation of Lemma~\ref{lem:rooted-threshold}. Let
		\(\lambda_1=\lambda_1(\Gamma)\),
		\(L=c-1,\)
		\(D=\sum_i d_i,\) and \(\eps=\lambda_1-L.\)
		Then
		\[\eps\lambda_1(\lambda_1+1)
		\le \sum_i d_i^2+\eps D+\frac{h(L-h)}{\lambda_1+1}.\]
	\end{lemma}
	
	\begin{proof}
		Since $C$ induces an all-positive $K_c$ in $\Gamma$,
		interlacing gives $\lambda_1\geq \lambda_1 (K_c)=c-1=L$. Thus
		$\varepsilon=\lambda_1-L\geq0$.
		We may assume that $x_{y_1}=1$. Then $x_y\leq1$ for every
		$y\in C$. Since $y_1\ne a$ and $y_1$ is adjacent to every
		vertex of $I$, the eigenvalue equation at $y_1$ gives
		\(
		\sum_{j=1}^{c}x_{y_j}
		=\lambda_1+1-\sum_{i=1}^{t}x_{z_i}.
		\)
		Summing the eigenvalue equations over $C$, we obtain
		\(
		\varepsilon\sum_{j=1}^{c}x_{y_j}
		=\sum_{i=1}^{t}d_i x_{z_i}+x_w,
		\)
		where $x_w$ comes from the positive edge $wa$.
		Combining these two equations gives
		\(
		\varepsilon(\lambda_1+1)
		=\sum_{i=1}^{t}(d_i+\varepsilon)x_{z_i}+x_w.
		\)
		For $1\leq i<t$, the eigenvalue equation at $z_i$ gives
		\(
		\lambda_1 x_{z_i}
		=\sum_{j=1}^{d_i}x_{y_j}\leq d_i.
		\)
		Since $z_t=b$ and $d_t=h$, the eigenvalue equation at $b$ gives
		\(
		\lambda_1 x_b+x_w=\sum_{j=1}^{h}x_{y_j}\leq h.
		\)
		Since $x_w\leq x_b$, it follows that
		$(\lambda_1+1)x_w\leq h$. Using these bounds and
		$\varepsilon\geq0$, we obtain
		\begin{align*}
			\varepsilon\lambda_1(\lambda_1+1)
			&=\sum_{i=1}^{t-1}(d_i+\varepsilon)\lambda_1 x_{z_i}
			+(h+\varepsilon)\lambda_1 x_b+\lambda_1 x_w\\
			&\leq\sum_{i=1}^{t-1}(d_i+\varepsilon)d_i
			+(h+\varepsilon)(h-x_w)+\lambda_1 x_w\\
			&=\sum_{i=1}^{t}d_i^2+\varepsilon D+(L-h)x_w\\
			&\leq\sum_{i=1}^{t}d_i^2+\varepsilon D
			+\frac{h(L-h)}{\lambda_1+1},
		\end{align*}
		where the last inequality follows from $h\leq L$ and
		$x_w\leq h/(\lambda_1+1)$. This proves Lemma \ref{lem:split-ear}.
	\end{proof}

	\begin{proof}[Proof of Theorem~\ref{t1.1}]
		Let \(\dot G=(G,\sigma)\) be an unbalanced \(C_3^-\)-free
		signed graph with \(m\) edges and maximum index among all
		such graphs. By Lemma~\ref{lem:one-component}, \(G\) is connected.
		By Lemma~\ref{lem:one-negative-edge}, there is a signed graph
		\(\dot G^\ast\) switching equivalent to \(\dot G\) with a
		unique negative edge \(wb\) and a nonnegative unit eigenvector
		\(\mathbf{x}\) corresponding to \(\lambda_1(\dot G^\ast)\).
		Since switching preserves the index and the sign of every
		cycle, \(\dot G^\ast\) is also extremal. Denote
		\(\lambda_1=\lambda_1(\dot G^\ast)\) and label \(w\) and \(b\)
		so that \(x_w\leq x_b\).
		
		We first prove that \(w\) has exactly one positive neighbor.
		Let \(H=G-wb\). Since \(\dot G^\ast\) is unbalanced, its unique
		negative edge \(wb\) lies on a negative cycle. Thus \(H\)
		is connected. Moreover, since \(\dot G^\ast\)
		is \(C_3^-\)-free, we have
		\(N_H(w)\cap N_H(b)=\varnothing\). 
		If \(x_b=0\), then \(x_w=0\), so
		\( A(H)\mathbf{x}
		= A(\dot G^\ast)\mathbf{x}
		= \lambda_1\mathbf{x}\). Since \(H\) is
		connected, the Perron--Frobenius theorem gives
		\(\mathbf{x}>\mathbf{0}\), a contradiction. Hence \(x_b>0\).
		Let \(P\) be a \(w\)--\(b\) path in \(H\) and \(a\) be
		the first vertex after \(w\) on \(P\). Suppose that
		\(u\in N_H(w)\setminus\{a\}\). Then \(bu\notin E(G)\).
		Let \(\dot G'\) be obtained from \(\dot G^\ast\) by replacing
		the positive edge \(wu\) with the positive edge \(bu\).
		Clearly, \(\dot G'\) still has \(m\) edges and the path \(P\) together with
		\(wb\) still forms a negative cycle. Moreover, \(\dot G'\) is \(C_3^-\)-free.
		By extremality and the Rayleigh principle,
		\(
		\lambda_1
		\geq \lambda_1(\dot G')
		\geq \mathbf{x}^{\mathsf T}A(\dot G')\mathbf{x}
		= \mathbf{x}^{\mathsf T}A(\dot G^\ast)\mathbf{x}
		-2x_wx_u+2x_bx_u
		= \lambda_1+2x_u(x_b-x_w)
		\geq \lambda_1.
		\)
		Thus equality holds throughout and \(\mathbf{x}\) is also
		an eigenvector of \(A(\dot G')\) corresponding to \(\lambda_1\).
		Comparing the eigenvalue equations at \(w\) in \(\dot G^\ast\)
		and \(\dot G'\) gives \(x_u=0\). Since all edges incident with
		\(u\) in \(\dot G'\) are positive, the eigenvalue equation at
		\(u\) gives
		\(
		0=\lambda_1 x_u
		=\sum_{v\in N_{G'}(u)}x_v
		\geq x_b>0,
		\)
		a contradiction. Hence \(N_H(w)=\{a\}\).

		Let \(J=H-w=G-w\). Since \(w\) is a leaf of \(H\), the graph
		\(J\) is connected. Deleting \(w\) removes exactly \(wa\) and
		\(wb\), so
		\(
		p:=e(J)=m-2=\binom{k}{2}+q-1.
		\)
		Moreover, since \(\dot G^\ast\) is \(C_3^-\)-free, we have \(ab\notin E(J)\). Thus
		\(\dot G^\ast=T(J;a,b)\).
		Since the candidate graphs \(\dot\Sigma_{k,q}\) for \(2\leq q\leq k\)
		and \(\dot\Omega_{k+1}\) for \(q=k+1\) each contain an induced \(C_4^-\),
		extremality and interlacing give
		\(\lambda_1\geq\lambda_1(C_4^-)=\sqrt{2}>1\).
		Since \(x_w\leq x_b\) and \(x_b>0\),
		by Lemma~\ref{lem:rooted-threshold}, we may
		write \(V(J)=C\mathbin{\dot\cup}I\), where
		\(C=\{y_1,\ldots,y_c\}\) is a maximum clique containing \(a\)
		and \(I=\{z_1,\ldots,z_t=b\}\) is an independent set. 
		Moreover,
		we have
		\(N_J(z_i)=\{y_1,\ldots,y_{d_i}\}\) for \(1\leq i\leq t\)
		and
		\(c-1\geq d_1\geq\cdots\geq d_t=h\geq1.\)

		Let \(L=c-1\) and \(D=\sum_{i=1}^t d_i=p-\binom{L+1}{2}\).
		Denote \(D=zL+r_0\), where \(0\leq r_0<L\). Since \(C\) induces
		an all-positive clique, \(\lambda_1\geq L\). Let
		\(\alpha=\sum_{i=1}^t d_i^2+\frac{h(L-h)}{\lambda_1+1}\)
		and
		\(F_D(X)=(X-L)\bigl[X(X+1)-D\bigr]\).
		If \(D\leq L+1\), then \(F_D(X)\) is strictly increasing on
		\((L,\infty)\), since both factors are positive and strictly
		increasing for \(X>L\).
		By Lemma \ref{lem:split-ear}, \(F_D(\lambda_1)\leq\alpha\).
		Since \(d_i\leq d_1\) for every \(1\leq i\leq t\) and
		\(D=\sum_{i=1}^t d_i\), we have
		\(\sum_{i=1}^t d_i^2\leq d_1D\).
		Moreover, since \(\lambda_1\geq L\) and \(h>0\), we have
		\(\frac{h(L-h)}{\lambda_1+1}<h\).
		When \(t\geq2\), we have \(h=d_t\leq D-d_1\), and hence
		\begin{align}\label{a3.5}
			F_D(\lambda_1)\leq\alpha<d_1D+h\leq d_1(D-1)+D.
		\end{align}
		
		We claim that
		\(
		\sum_{i=1}^{t}d_i^2\leq zL^2+r_0^2.
		\)
		Indeed, transferring weight from a smaller entry to a larger one
		increases the sum of squares until one of the two entries reaches
		\(0\) or \(L\). Repeating this operation leaves at most one entry strictly between
		\(0\) and \(L\). Since the total sum remains \(D=zL+r_0\), the
		resulting sequence has \(z\) entries equal to \(L\), one entry
		equal to \(r_0\) if \(r_0>0\), and zeros elsewhere. Its sum of
		squares is therefore \(zL^2+r_0^2\). Thus
		\(
		\sum_{i=1}^t d_i(L-d_i)
		=LD-\sum_{i=1}^t d_i^2
		\geq LD-(zL^2+r_0^2)
		=r_0(L-r_0).
		\)
		Since \(0\leq d_i\leq L\) for each \(i\), we have
		\(
		\sum_{i=1}^{t}d_i(L-d_i)
		\geq d_t(L-d_t)=h(L-h),
		\)
		where the equality follows from \(d_t=h\). Moreover, since
		\(D=\sum_{i=1}^{t}d_i\), we have
		\(
		\sum_{i=1}^{t}d_i^2
		=LD-\sum_{i=1}^{t}d_i(L-d_i).
		\)
		By the definition of \(\alpha\) and the preceding estimates, we have
		\[
		\alpha
		=LD-\sum_{i=1}^{t}d_i(L-d_i)
		+\frac{h(L-h)}{\lambda_1+1}
		\leq LD-\frac{\lambda_1}{\lambda_1+1}
		\sum_{i=1}^{t}d_i(L-d_i)
		\leq LD-\frac{\lambda_1 r_0(L-r_0)}{\lambda_1+1}.\]
		Combining this with \(F_D(\lambda_1)\leq\alpha\) gives
		\begin{align}\label{a3.4}
			g(\lambda_1)\leq0,
		\end{align}
		where
		\(
		g(X)=(X-L)(X+1)-D+\frac{r_0(L-r_0)}{X+1}.
		\)
		For \(X\geq L\), we have \(2X-L+1\geq L+1\) and
		\(X+1\geq L+1\). Since \(r_0(L-r_0)\leq \frac{L^2}{4}\), we obtain
		\(
		g'(X)
		=2X-L+1-\frac{r_0(L-r_0)}{(X+1)^2}
		\geq L+1-\frac{L^2}{4(L+1)^2}>0.
		\)
		Thus \(g(X)\) is strictly increasing on \([L,\infty)\).

				\medskip
		\noindent\textbf{Case 1.} \textit{\(2\leq q\leq k\).}
		
		Let \(s=q-1\) and \(R=\lambda_1(R_p)\). By extremality and
		Lemma~\ref{lem:candidate-lower}\textup{(i)}, \(\lambda_1
		\geq\lambda_1(\dot\Sigma_{k,q})
		\geq R\).
		Since \(p=\binom{k}{2}+s\leq \binom{k}{2}+k-1
		<\binom{k+1}{2}\), we have \(c\leq k\).
		Let \(\eta=R-k+1>0\).
		
		Let \(v\) be the added vertex of \(R_p\) and
		\(\widetilde{\mathbf{x}}\) be a positive eigenvector of \(A(R_p)\)
		corresponding to \(R\).
		Let \(S=\sum_{u\in V(K_k)}\widetilde{x}_u\).
		Summing the eigenvalue equations over the vertices of \(K_k\)
		gives
		\(
		RS
		=\sum_{u\in V(K_k)}(S-\widetilde{x}_u)+s\widetilde{x}_v
		=(k-1)S+s\widetilde{x}_v.
		\)
		Thus \(\eta S=s\widetilde{x}_v\).
		For each \(u\in N_{R_p}(v)\), the eigenvalue equation at \(u\)
		gives 
		\((R+1)\widetilde{x}_u=S+\widetilde{x}_v\).
		Summing over \(N_{R_p}(v)\) and using the eigenvalue equation at \(v\),
		namely \(\sum_{u\in N_{R_p}(v)}\widetilde{x}_u
		=R\widetilde{x}_v\), we obtain
		\(
		R(R+1)\widetilde{x}_v
		=s(S+\widetilde{x}_v).
		\)
		Multiplying by \(\eta\) and substituting \(\eta S=s\widetilde{x}_v\),
		we get
		\(
		\eta R(R+1)\widetilde{x}_v
		=s(\eta S+\eta\widetilde{x}_v)
		=s(s+\eta)\widetilde{x}_v.
		\)
		Since \(\widetilde{x}_v>0\), it follows that
		\begin{align}\label{a3.2}
			\eta R(R+1)=s(s+\eta).
		\end{align}
		
		We claim that \(c=k\). Suppose to the contrary that \(c<k\). Let \(n=k-1-L=k-c\geq1\),
		\(M=\binom{n+1}{2}\), and \(\theta=\frac{s+\eta}{R}\).
		Since \(s\leq k-1\) and \(R=k-1+\eta\), we have
		\(\theta\leq1\). Using \(k=L+n+1\), we obtain
		\(
		D
		=\binom{k}{2}-\binom{L+1}{2}+s
		=nL+M+s.
		\)
		Dividing \eqref{a3.2} by \(R\) gives
		\(\eta(R+1)=s\theta\). Hence
		\begin{align}\label{a3.3}
			s-\eta R=s+\eta-s\theta=(R-s)\theta.
		\end{align}
		By the definition of \(g\) and the identities
		\(R-L=n+\eta\) and \(D=nL+M+s\), we have
		\begin{align}
			g(R)
			&=(n+\eta)(R+1)-(nL+M+s)
			+\frac{r_0(L-r_0)}{R+1}\notag\\
			&=M+n\eta+s(\theta-1)+\frac{r_0(L-r_0)}{R+1},
			\label{eq:g-lambda1-prime}
		\end{align}
		where the last equality follows from
		\(n(R+1-L)=2M+n\eta\) and
		\(\eta(R+1)=s\theta\).
		
		Since \(R\theta=s+\eta\) and \(R=L+n+\eta\),
		we have \(R(\theta-1)=s-L-n\). Multiplying both sides of \eqref{eq:g-lambda1-prime} by \(R\)
		and using the identity \(\eta R=s-(R-s)\theta\),
		which follows from \eqref{a3.3}, we obtain
		\begin{align*}
			Rg(R)
			=(R-s)(M-n\theta)+s(M+s-L)
			+\frac{Rr_0(L-r_0)}{R+1}.
		\end{align*}
		
		We next prove that \(g(R)>0\).
		Since \(R>k-1\geq s\), \(M\geq n\), and \(\theta\leq1\),
		the first term on the right-hand side is nonnegative.
		The third term is also nonnegative,
		since \(0\leq r_0<L\).
		If \(M+s>L\), then the second term is
		positive. If \(M+s=L\), then \(s=L-M<L<k-1\), and hence
		\(\theta=\frac{s+\eta}{R}=\frac{s+\eta}{k-1+\eta}<1\). Thus
		\(M-n\theta\geq n(1-\theta)>0\), so the first term is positive.
		Finally, if \(M+s<L\), then \(D=nL+M+s\) gives \(r_0=M+s\).
		The last two terms therefore sum to
		\(
		\frac{(L-r_0)(RM-s)}{R+1}>0,
		\)
		where the strict inequality follows from \(r_0<L\) and
		\(RM\geq R>s\).
		In each case, \(Rg(R)>0\), so
		\(g(R)>0\).
		Since \(g\) is increasing for \(X\geq L\)
		and \(\lambda_1\geq R>L\), we have
		\(g(\lambda_1)\geq g(R)>0\), contradicting
		\eqref{a3.4}. Hence \(c=k\) and
		\(D=p-\binom{k}{2}=s\leq k-1=L\).

		If \(t\geq2\), substituting \(D=s\) into \eqref{a3.5} and using \(d_1\leq s\),
		we obtain
		\(
		F_s(\lambda_1)\leq\alpha<d_1(s-1)+s\leq s^2.
		\)
		On the other hand, since \(s\leq L\), \(F_s(X)\) is strictly
		increasing for \(X>L\). Thus \(\lambda_1\geq R>L\)
		gives \(F_s(\lambda_1)\geq F_s(R)\).
		Since \(L=k-1\), we have \(R-L=\eta\).
		By \eqref{a3.2}, we obtain
		\(
		F_s(R)
		=\eta R(R+1)-\eta s
		=s^2.
		\)
		Hence \(F_s(\lambda_1)\geq s^2\), a contradiction. 
		Since \(b\in I\), we must have \(t=1\) and \(I=\{b\}\).
		Thus \(d_J(b)=D=s=q-1\).
		Since \(a\in C\setminus N_J(b)\), we obtain
		\(\dot G^\ast\cong\dot\Sigma_{k,q}\).

		\medskip
		\noindent\textbf{Case 2.} \textit{\(q=k+1\).}
		
		Let \(d=k+1\) and \(\tau=\lambda_1(\dot\Omega_d)\).
		Then \(p=\binom{d}{2}\) and \(\lambda_1\geq\tau\).
		Since \(D\geq1\), we have \(\binom{c}{2}<p=\binom{d}{2}\).
		Thus \(c\leq d-1\), and hence \(2\leq c\leq d-1\).
		We first claim that \(c=d-1\).
		Suppose that \(c\leq d-2\). Then \(d\geq4\) and \(L\leq d-3\).
		For \(1\leq\ell\leq d-3\), let
		\(
		\phi(\ell)=(\tau-\ell)(\tau+1)
		-\left(\binom{d}{2}-\binom{\ell+1}{2}\right).
		\)
		For \(1\leq\ell<d-3\), we have \(\tau>d-2>\ell\), so
		\(\phi(\ell+1)-\phi(\ell)=\ell-\tau<0\).
		Since \(L\leq d-3\), it follows that \(\phi(L)\geq\phi(d-3)\).
		By
		Lemma~\ref{lem:candidate-lower}\textup{(ii)}, we obtain
		\(
		\tau(\tau-d+3)=2(d-2)+\frac{2}{\tau}.
		\)
		Since \(\binom{d}{2}-\binom{d-2}{2}=2d-3\), the definition of
		\(\phi\) gives
		\(\phi(d-3)=\tau-(d-2)+\frac{2}{\tau}>0\),
		where the inequality follows from \(\tau>d-2>0\).
		Substituting \(D=\binom{d}{2}-\binom{L+1}{2}\) into the
		definition of \(g\), we obtain
		\(
		g(\tau)
		=\phi(L)+\frac{r_0(L-r_0)}{\tau+1}
		\geq\phi(L)\geq\phi(d-3)>0,
		\)
		where the first inequality follows from \(0\leq r_0<L\).
		Since \(\lambda_1\geq\tau>d-2\geq c>L\) and \(g\) is increasing
		for \(X\geq L\), we have \(g(\lambda_1)\geq g(\tau)>0\),
		contradicting \eqref{a3.4}.
		Thus \(c= d-1\), \(L=c-1=d-2\), and \(D=p-\binom{c}{2}
		=\binom{d}{2}-\binom{d-1}{2}
		=L+1.\)

		If \((d_1,\ldots,d_t)\ne(L,1)\), then \(t\geq2\) and
		\(d_1\leq L-1\), so \(L\geq2\). Substituting \(D=L+1\) into \eqref{a3.5} and using \(d_1\leq L-1\),
		we obtain
		\(
		F_{L+1}(\lambda_1)\leq\alpha
		<d_1L+L+1\leq L^2+1.
		\)
		On the other hand, by the definition of \(F_{L+1}\) and the cubic
		equation for \(\tau\) in Lemma~\ref{lem:candidate-lower}\textup{(ii)},
		we have
		\(
		F_{L+1}(\tau)=L(L+1)+2-\tau>L^2+1,
		\)
		where the last inequality follows from \(\tau<d-1=L+1\).
		Since \(F_{L+1}(X)\) is increasing for \(X>L\) and
		\(\lambda_1\geq\tau>L\), we have
		\(
		F_{L+1}(\lambda_1)\geq F_{L+1}(\tau)>L^2+1,
		\)
		a contradiction. Therefore \((d_1,\ldots,d_t)=(L,1)\).

		Thus \(I=\{z_1,b\}\),
		\(N_J(z_1)=\{y_1,\ldots,y_L\}\), and \(N_J(b)=\{y_1\}\).
		Hence \(J\cong L_d\). Since \(ab\notin E(J)\), we have
		\(a\ne y_1\). If \(d=3\), this gives \(a=y_2\) and
		\(\dot G^\ast\cong\dot\Omega_3\).
		Now let \(d\geq4\). If \(a=y_{L+1}\), then \(by_2\notin E(J)\)
		and the argument in the proof of Lemma~\ref{lem:rooted-threshold}
		gives \(x_a\geq x_{y_2}\). 
		The eigenvalue equations at \(z_1\) and \(b\) give
		\(
		\lambda_1x_{z_1}\geq x_{y_1}=\lambda_1x_b+x_w.
		\)
		Since \(\lambda_1>0\), it follows that
		\(x_{z_1}\geq x_b+\frac{x_w}{\lambda_1}\).
		If \(x_b>x_w\), then \(x_{z_1}\geq x_b>x_w\).
		If \(x_b=x_w\), then \(x_w=x_b>0\), so
		\(x_{z_1}\geq x_w+\frac{x_w}{\lambda_1}>x_w\).
		Thus \(x_{z_1}>x_w\) in either case.
		By the eigenvalue equations at \(a\) and \(y_2\), we have
		\(
		(\lambda_1+1)(x_{y_2}-x_a)=x_{z_1}-x_w>0.
		\)
		Hence \(x_{y_2}>x_a\), contradicting \(x_a\geq x_{y_2}\).
		Therefore \(a\in\{y_2,\ldots,y_L\}\) and
		\(\dot G^\ast\cong\dot\Omega_d\).
		
		In both cases, every extremal graph is switching isomorphic to the
		corresponding candidate. Conversely, switching isomorphism preserves
		the index, so every graph switching isomorphic to that candidate
		attains equality.
	\end{proof}
	
	\subsection{A complete characterization of the adjacency spectral radius extremal graphs}
	
	Since isolated vertices do not affect the spectral radius, we may
	assume throughout that \(\dot G\) has no isolated vertices.
	Let \(\Lambda_m\) denote the maximum value of \(\lambda_1(\dot G)\)
	over all unbalanced \(C_3^-\)-free signed graphs \(\dot G\)
	with exactly \(m\) edges, as determined in Theorem~\ref{t1.1}.
	\begin{lemma}\label{c4.1}
		For every \(m\geq7\), we have \(\Lambda_m>\sqrt{m-1}\).
	\end{lemma}
	
	\begin{proof}
		For \(7\leq m\leq12\), Theorem~\ref{t1.1} gives the graphs
		attaining the maximum index \(\Lambda_m\). We compute their signed
		adjacency characteristic polynomials directly and evaluate the
		monic factors \(f_m(t)\) at \(t=\sqrt{m-1}\). The results are
		shown in the following table.
		Since \(f_m(t)\to+\infty\) as \(t\to+\infty\) and
		\(f_m(\sqrt{m-1})<0\), we have \(\Lambda_m>\sqrt{m-1}\).
		\begin{center}
			\begin{tabular}{@{}ccll@{}}
				\toprule
				\(m\) & Candidate & \(f_m(t)\) & \(f_m(\sqrt{m-1})\) \\
				\midrule
				\(7\)  & \(\dot\Sigma_{3,3}\) & \(t^2-t-4\) & \(2-\sqrt{6}<0\) \\
				\(8\)  & \(\dot\Omega_4\)     & \(t^2-2t-2\) & \(5-2\sqrt{7}<0\) \\
				\(9\)  & \(\dot\Sigma_{4,2}\) & \(t^3-t^2-6t-2\) & \(4\sqrt{2}-10<0\) \\
				\(10\) & \(\dot\Sigma_{4,3}\) & \(t^5-t^4-9t^3-t^2+15t+7\) & \(-38<0\) \\
				\(11\) & \(\dot\Sigma_{4,4}\) & \(t^2-2t-6\) & \(4-2\sqrt{10}<0\) \\
				\(12\) & \(\dot\Omega_5\)     & \(t^3-2t^2-6t-2\) & \(5\sqrt{11}-24<0\) \\
				\bottomrule
			\end{tabular}
		\end{center}

		Now suppose that \(m\geq13\). By the definition of \(k\), we have
		\(k\geq5\). For \(\dot\Sigma_{k,q}\), let \(J\) be the all-positive
		induced subgraph obtained by deleting \(u\). It consists of \(K_k\)
		and a vertex \(v\) with \(q-1\geq1\) neighbors in \(K_k\). For
		\(\dot\Omega_{k+1}\), let \(J=K_{k+1}-xy\) be the all-positive
		induced subgraph obtained by deleting \(u\) and \(v\). Deleting
		\(x\) from \(J\) leaves a \(K_k\).
		In both cases, \(J\) is a connected all-positive induced subgraph
		of the corresponding candidate and properly contains \(K_k\).
		By interlacing and the strict monotonicity of the spectral radius
		for connected unsigned graphs, we obtain
		\(
		\Lambda_m\geq\rho(J)>\rho(K_k)=k-1.
		\)
		Moreover,
		\(
		m-1\leq\binom{k+1}{2}+1
		=(k-1)^2-\frac{k(k-5)}{2}
		\leq(k-1)^2.
		\)
		Hence \(\Lambda_m>k-1\geq\sqrt{m-1}\).
	\end{proof}
	\begin{lemma}\label{c4.2}
		For every \(m\geq7\), we have
		\(-\lambda_n(\dot G)\leq\sqrt{m-1}\).
	\end{lemma}
	
	\begin{proof}
		Let \(\mu=-\lambda_n(\dot G)\) and \(-\dot G\) be
		obtained from \(\dot G\) by reversing the sign of every edge.
		Since \(A(-\dot G)=-A(\dot G)\), we have
		\(\lambda_1(-\dot G)=\mu\). By Lemma~\ref{l1}, there exists
		a signed graph \(\widetilde G\) switching equivalent to
		\(-\dot G\) such that \(A(\widetilde G)\) has a nonnegative
		unit eigenvector \(\mathbf{x}\) corresponding to \(\mu\).
		
		Let \(H\) be the spanning subgraph of \(G\) consisting of
		the positive edges of \(\widetilde G\). Since \(\dot G\) is
		\(C_3^-\)-free, every triangle in \(-\dot G\) is negative.
		Switching preserves the sign of every cycle, so every
		triangle in \(\widetilde G\) is negative. Thus \(H\) is \(C_3\)-free.
		Since \(\mathbf{x}\geq\mathbf{0}\), the Rayleigh principle gives
		\(
		\mu=\mathbf{x}^{\mathsf T}A(\widetilde G)\mathbf{x}
		\leq\mathbf{x}^{\mathsf T}A(H)\mathbf{x}
		\leq\rho(H).
		\)
		
		For each \(v\in V(H)\), the \(v\)-th row sum of \(A(H)^2\)
		is \(\sum_{u\in N_H(v)}d_H(u)\). Since \(N_H(v)\) is an
		independent set, this sum counts each edge with an endpoint
		in \(N_H(v)\) exactly once. Hence the row sum is
		\(
		e(H)-e\bigl(H[V(H)\setminus N_H(v)]\bigr).
		\)
		If \(\widetilde G\) has a negative edge, then \(e(H)\leq m-1\).
		Otherwise, \(e(H)=m\) and \(\dot G\) is switching equivalent
		to the all-negative signing of \(H\). Since \(\dot G\) is
		unbalanced, \(H\) contains an odd cycle and hence is not
		bipartite. Thus \(V(H)\setminus N_H(v)\) contains an edge
		for every \(v\in V(H)\); otherwise, it and \(N_H(v)\) would
		form a bipartition of \(H\).
		
		In either case, every row sum of \(A(H)^2\) is at most
		\(m-1\). 
		Since \(A(H)^2\) is nonnegative, its spectral radius
		is at most its maximum row sum. 
		Therefore
		\(\rho(H)^2\leq m-1\), and hence
		\(\mu\leq\rho(H)\leq\sqrt{m-1}\).
	\end{proof}
	
	\begin{lemma}\label{c4.3}
		If \(m=6\), then \(-\lambda_n(\dot G)<\Lambda_6\).
	\end{lemma}
	
	\begin{proof}
		When \(m=6\), we have \(k=3\) and \(q=2\). Hence
		Theorem~\ref{t1.1} gives
		\(\Lambda_6=\lambda_1(\dot\Sigma_{3,2})\).
		Suppose first that \(G\) is bipartite. By the symmetry of the
		signed adjacency spectrum about zero, we have
		\(
		\sum_{\lambda_i(\dot G)>0}\lambda_i^2(\dot G)
		=\frac12\operatorname{tr}\bigl(A(\dot G)^2\bigr)
		=m=6.
		\)
		If \(\dot G\) contains a \(C_4^-\), then interlacing gives
		\(
		\lambda_2(\dot G)\geq\lambda_2(C_4^-)=\sqrt{2}.
		\)
		Hence
		\(
		\rho(\dot G)^2=\lambda_1(\dot G)^2
		\leq 6-\lambda_2(\dot G)^2
		\leq 4.
		\)
		Thus \(\rho(\dot G)\leq 2\).
		If \(\dot G\) is \(C_4^-\)-free, then it contains a \(C_6^-\). Direct computation gives
		\(\rho(\dot G)=\rho(C_6^-)=\sqrt{3}\).

		We claim that if \(G\) is nonbipartite, then it contains a \(C_5\).
		Suppose that \(G\) is nonbipartite and contains no \(C_5\).
		Then \(G\) contains a triangle. Since \(\dot G\) is unbalanced
		and \(C_3^-\)-free, it contains a \(C_4^-\). 
		Since \(m=6\), the triangle shares an edge with the \(C_4^-\).
		A chord of the \(C_4^-\) would create a negative triangle, so
		the third vertex of the triangle lies outside the \(C_4^-\).
		Removing the common edge from their union yields a \(C_5\),
		a contradiction.
		
		Thus \(\dot G\) contains a \(C_5^-\). Let \(e\) be the remaining edge. There are three cases.
		If \(e\) forms a separate component, then \(
		\rho(\dot G)=\max\{\rho(C_5^-),1\}=\rho(C_5^-)=2.
		\)
		If \(e\) is a pendant edge of the \(C_5^-\), we may switch
		\(\dot G\) so that all edges are negative. Direct computation
		gives the eigenvalues \(-1\), \(\frac{1\pm\sqrt{5}}{2}\), and
		the three roots of \(t^3-4t+1=0\). Hence \(\rho(\dot G)=\beta\),
		where \(\beta\) is the unique root of \(t^3-4t-1=0\) in
		\((2,\infty)\).
		If \(e\) is a chord, then it forms a \(C_3^+\) and a
		\(C_4^-\). Thus \(\dot G\) is switching isomorphic
		to \(\dot\Sigma_{3,2}\). Direct computation gives the
		eigenvalues \(\pm\sqrt{2}\) and the three roots of
		\(t^3-4t-2=0\). The only eigenvalue outside \([-2,2]\) is
		the unique root \(\alpha\) in \((2,\infty)\), so
		\(\rho(\dot G)=\alpha=\Lambda_6\).
		Since \(\beta^3-4\beta-2=-1<0\) and \(t^3-4t-2\) is strictly
		increasing on \([2,\infty)\), we have \(2<\beta<\alpha\).
		Together with the bipartite case, this gives
		\(
		-\lambda_n(\dot G)\leq\beta<\Lambda_6.
		\)
	\end{proof}
	
	\begin{proof}[Proof of Theorem \ref{t2}]
		Since \(\dot G\) is unbalanced and \(C_3^-\)-free, it contains
		a negative cycle of length at least four.
		If \(m=4\), then \(\dot G\) is a \(C_4^-\). Thus \(\dot G\) is switching isomorphic
		to \(\dot\Sigma_{2,2}\), the extremal graph given by
		Theorem~\ref{t1.1} for \(m=4\). Its eigenvalues are
		\(\pm\sqrt{2}\), each with multiplicity two, so
		\(\rho(\dot G)=\sqrt{2}\).
		Suppose that \(m=5\). If \(\dot G\) contains a \(C_5^-\), then it is
		switching isomorphic to \((C_5,-)\). The eigenvalues of \((C_5,-)\)
		are \(-2\) and \(\frac{1\pm\sqrt{5}}{2}\), the latter two each with
		multiplicity two, so \(\rho(\dot G)=2\). Otherwise, \(\dot G\)
		contains a \(C_4^-\), and the remaining edge cannot be a chord,
		since it would create a \(C_3^-\). If it is a pendant edge, the
		resulting graph is switching isomorphic to \(\dot\Omega_3\), with
		eigenvalues \(0,\pm\sqrt{2},\pm\sqrt{3}\). Thus
		\(\rho(\dot G)=\sqrt{3}\). If the edge forms a separate component,
		then \(\rho(\dot G)=\sqrt{2}\). Hence \(\rho(\dot G)\leq2\), with
		equality precisely when \(\dot G\) is switching isomorphic to
		\((C_5,-)\).
		
		It remains to consider \(m\geq6\). By Theorem~\ref{t1.1},
		we have \(\lambda_1(\dot G)\leq\Lambda_m\). If \(m=6\),
		by Lemma~\ref{c4.3},
		we have
		\(-\lambda_n(\dot G)<\Lambda_6\). If
		\(m\geq7\), then Lemmas~\ref{c4.1} and~\ref{c4.2} give
		\(
		-\lambda_n(\dot G)\leq\sqrt{m-1}<\Lambda_m.
		\)
		Thus \(-\lambda_n(\dot G)<\Lambda_m\) for every \(m\geq6\), and hence
		\(
		\rho(\dot G)
		=\max\{\lambda_1(\dot G),-\lambda_n(\dot G)\}
		\leq\Lambda_m,
		\)
		with equality if and only if \(\lambda_1(\dot G)=\Lambda_m\).
		By Theorem~\ref{t1.1}, this holds if and only if \(\dot G\)
		is switching isomorphic to \(\dot\Sigma_{k,q}\) when \(2\leq q\leq k\),
		or to \(\dot\Omega_{k+1}\) when \(q=k+1\).
	\end{proof}
	
\section{Spectral extremal results for unbalanced
	\(\mathcal K_r^-\)-free signed graphs with \(r\geq4\)}

\subsection{A complete characterization of the index extremal graphs}

\begin{proof}[Proof of Theorem~\ref{thm:Kr-index}.]
	
Let \(r\geq4\), \(m\geq3\), and 
\(p=m-2=\binom{k}{2}+s\), where \(0\le s\le k-1\) and \(k\geq2\).
We first consider the range \(m\geq7\) and \(1\leq s\leq M\), where \(M=\min\{k-2,r-4\}\).
If \(s\leq k-3\), then
\(
m=\binom{k}{2}+(s+1)+1.
\)
Applying Theorem \ref{w}(3) with
\(d=k\) and \(t=s+1\), we have every extremal graph among all unbalanced signed
graphs with \(m\) edges is switching isomorphic to
\(\dot\Theta_{k,s}\) or \(\dot B_{k,s+1}\).
If \(s=k-2\), then \(m=\binom{k+1}{2}\), so
Theorem \ref{w}(1) gives the unique
extremal graph \(\dot B_{k,k-1}\).
Since \(s\leq M\leq r-4\), all these extremal graphs are
\(\mathcal K_r^-\)-free. Thus the maximum index among
unbalanced \(\mathcal K_r^-\)-free signed graphs with \(m\)
edges equals that among all unbalanced signed graphs with
\(m\) edges. The comparison of the two
candidates will be given at the end of this proof.
	
	Let \(\dot G=(G,\sigma)\) be an unbalanced
	\(\mathcal K_r^-\)-free signed graph with \(m\) edges
	and maximum index among all such graphs.
	We may assume that \(\dot G\) has no isolated vertices.
	By Lemma~\ref{lem:one-component}, \(\dot G\) is connected.
	By Lemma~\ref{lem:one-negative-edge}, there exists a signed
	graph \(\dot G^\ast\) switching equivalent to \(\dot G\)
	with exactly one negative edge \(ub\), such that
	\(A(\dot G^\ast)\) has a nonnegative unit eigenvector
	\(\mathbf{x}\) corresponding to
	\(\lambda=\lambda_1(\dot G^\ast)\).
	Interchanging \(u\) and \(b\) if necessary, we may assume
	that \(x_u\leq x_b\).
	Since switching preserves the index and the sign of every
	cycle, \(\dot G^\ast\) is also extremal.
	We claim that \(x_b>0\). Otherwise, \(x_u=x_b=0\).
	Since \(ub\) is the only negative edge, we have
	\(A(G)\mathbf{x}=A(\dot G^\ast)\mathbf{x}
	=\lambda\mathbf{x}\).
	Since \(G\) is connected, the Perron--Frobenius theorem gives
	\(\mathbf{x}>\mathbf{0}\), a contradiction.
	Let \(P\) be the spanning
	subgraph of \(G\) consisting of the positive edges of
	\(\dot G^\ast\) and
	\(
	S=N_P(u)\cap N_P(b).
	\)

	\begin{lemma}\label{l5.1}
	 \(S\ne\varnothing\).
	\end{lemma}
	
\begin{proof}[Proof of Lemma~\ref{l5.1}]
Suppose to the contrary that \(S=\varnothing\).
Since \(ub\) lies on a negative cycle, choose a positive
\(u\)--\(b\) path \(ua\cdots yb\) in \(P\).
Suppose that \(z\in N_P(u)\setminus\{a\}\).
Then \(bz\notin E(G)\).
Let \(\dot G'\) be obtained from \(\dot G^\ast\) by
replacing the positive edge \(uz\) with the positive
edge \(bz\).
Clearly, \(\dot G'\) still has \(m\) edges, and the 
positive \(u\)--\(b\) path together with \(ub\) still
forms a negative cycle.
Moreover, \(\dot G'\) is \(\mathcal K_r^-\)-free, since \(N_{G'}(u)\cap N_{G'}(b)=\varnothing\).
By the Rayleigh principle,
we have
\(
	\Delta
	=\mathbf{x}^{\mathsf T}A(\dot G')\mathbf{x}
	-\mathbf{x}^{\mathsf T}A(\dot G^\ast)\mathbf{x}
	=2x_z(x_b-x_u)\geq0.
\)
We consider the following four cases.
\begin{enumerate}
	\renewcommand{\labelenumi}{(\arabic{enumi})}
	\item
	If \(x_b>x_u\) and \(x_z>0\), then \(\Delta>0\).
	By the Rayleigh principle,
	\(\lambda_1(\dot G')>\lambda_1(\dot G^\ast)\),
	contradicting extremality.
	
	\item
	If \(x_b=x_u\) and \(x_z>0\), then \(\Delta=0\), but
	\(
	\bigl((A(\dot G')-A(\dot G^\ast))
	\mathbf{x}\bigr)_u=-x_z<0.
	\)
	By Lemma~\ref{lem:strict-increase},
	\(\lambda_1(\dot G')>\lambda_1(\dot G^\ast)\),
	a contradiction.
	
	\item
	If \(x_b>x_u\) and \(x_z=0\), then \(\Delta=0\), but
	\(
	\bigl((A(\dot G')-A(\dot G^\ast))
	\mathbf{x}\bigr)_z=x_b-x_u>0.
	\)
	Again, Lemma~\ref{lem:strict-increase} gives a
	contradiction.
	
	\item
	If \(x_b=x_u\) and \(x_z=0\), then 
	the eigenvalue equation at \(z\) gives
	\(
	0=\lambda x_z
	=\sum_{v\in N_P(z)}x_v
	\geq x_u=x_b>0,
	\)
	a contradiction.
\end{enumerate}
Thus no vertex \(z\in N_P(u)\setminus\{a\}\) exists,
and hence \(N_P(u)=\{a\}\).
Since \(S=\varnothing\), we have \(ab\notin E(P)\)
and \(y\ne a\). Moreover, \(N_P(u)=\{a\}\) implies \(d_G(u)=2\).
If \(x_y>x_a\), let \(\dot G_1\) be obtained from
\(\dot G^\ast\) by replacing the positive edge \(ua\)
with \(uy\). Clearly, \(\dot G_1\) still has \(m\) edges
and contains the negative triangle \(uybu\).
Since \(d_G(u)=2\) and \(ub\) is the only negative
edge, \(\dot G_1\) is \(\mathcal K_r^-\)-free. We have
\(
\Delta_1
=\mathbf{x}^{\mathsf T}A(\dot G_1)\mathbf{x}
-\mathbf{x}^{\mathsf T}A(\dot G^\ast)\mathbf{x}
=2x_u(x_y-x_a)\geq0.
\)
If \(x_u>0\), then \(\Delta_1>0\), contradicting
extremality. If \(x_u=0\), then
\(
\bigl((A(\dot G_1)-A(\dot G^\ast))\mathbf{x}\bigr)_u
=x_y-x_a>0.
\)
Thus Lemma~\ref{lem:strict-increase} also gives a
contradiction.
If \(x_y\leq x_a\), let \(\dot G_2\) be obtained from
\(\dot G^\ast\) by replacing the positive edge \(by\)
with \(ba\).
A similar argument gives a contradiction.
Therefore \(S\ne\varnothing\).
\end{proof}
We claim that \(S=N_P(u)\).
Clearly, \(S\subseteq N_P(u)\).
Suppose to the contrary that there exists
\(z\in N_P(u)\setminus S\). Then \(bz\notin E(G)\).
Let \(\dot G'\) be obtained from \(\dot G^\ast\)
by replacing the positive edge \(uz\) with the positive
edge \(bz\).
Clearly, \(\dot G'\) still has \(m\) edges.
Moreover, \(N_{G'}(u)\cap N_{G'}(b)=S\), and all edges
within \(S\) are unchanged. Thus \(\dot G'\) is
\(\mathcal K_r^-\)-free.
By Lemma~\ref{l5.1}, \(S\ne\varnothing\).
For any \(w\in S\), the negative triangle \(ubwu\)
remains in \(\dot G'\), so \(\dot G'\) is unbalanced.
The same Rayleigh quotient argument as in the proof
of Lemma~\ref{l5.1} gives
\(\lambda_1(\dot G')>\lambda_1(\dot G^\ast)\),
contradicting the extremality of \(\dot G^\ast\).
Hence \(S=N_P(u)\).

	We next prove that \(S\) is a clique.
	We first show that \(x_w>x_u\) for every \(w\in S\).
	The eigenvalue equation at \(w\) gives \(x_w>0\),
	since \(x_b>0\).
	Suppose that \(x_w\leq x_u\). 
	Let \(\dot G''\) be obtained from \(\dot G^\ast\)
	by interchanging the signs of \(ub\) and \(uw\).
	Clearly, \(\dot G''\) still has \(m\) edges and contains
	the negative triangle \(ubwu\).
	We claim that \(\dot G''\) is \(\mathcal K_r^-\)-free.
	Suppose to the contrary that \(\dot G''\) contains an
	unbalanced \(K_r\), say \(Q\).
	Since \(uw\) is the only negative edge of \(\dot G''\),
	we have \(uw\in E(Q)\).
	If \(b\in V(Q)\), then the same vertices give an
	unbalanced \(K_r\) in \(\dot G^\ast\), with \(ub\)
	as its only negative edge.
	If \(b\notin V(Q)\), then
	\(V(Q)\setminus\{u\}\subseteq S\subseteq N_P(b)\).
	Replacing \(w\) with \(b\) gives an unbalanced \(K_r\)
	in \(\dot G^\ast\). Thus, \(\dot G''\) is \(\mathcal K_r^-\)-free.
	Since \(x_b\geq x_u\geq x_w>0\),
	the same argument as in the proof of Lemma~\ref{l5.1}
	gives a contradiction.
	Hence \(x_w>x_u\) for every \(w\in S\).
Suppose that \(z,w\in S\) are nonadjacent.
Replacing the positive edge \(uz\) with the positive
edge \(zw\) preserves the number of edges and the
negative triangle \(ubwu\).
The resulting graph is still \(\mathcal K_r^-\)-free.
Since \(x_z>x_u\ge0\) and \(x_w>x_u\), the same argument
as in the proof of Lemma~\ref{l5.1} gives a contradiction.
Therefore \(S\) is a clique.
	
	Let \(H=G-u=P-u\).
	Since \(G\) is connected and \(ub\) lies on a negative
	cycle, \(P\) is connected.
	Since \(N_P(u)=S\) is a clique, \(H\) is connected.
	Let \(t=|S|\). 
	Since \(S\) is a clique and \(\dot G^\ast\) is
	\(\mathcal K_r^-\)-free, we have
	\(1\leq t\leq r-3\) and \(d_G(u)=t+1\leq r-2.\)
	Choose \(a,v\in S\) such that
	\(x_a=\max_{w\in S}x_w\)
	and
	\(x_v=\min_{w\in S}x_w.\)
	For \(z\in V(H)\setminus(S\cup\{b\})\), replace the
	positive edge \(uv\) with \(uz\).
	If \(t\geq2\), a negative triangle through a vertex
	of \(S\setminus\{v\}\) remains.
	If \(t=1\), a positive \(z\)--\(b\) path in \(H\),
	together with \(uz\) and \(ub\), forms a negative cycle.
	If \(x_z>x_v\), the same argument as in the proof
	of Lemma~\ref{l5.1} gives a contradiction.
	Thus \(x_z\leq x_v\).
	Interchanging the signs of \(ub\) and \(uv\) preserves
	the negative triangle \(ubvu\).
	If \(x_b>x_v\), the same argument as in the proof
	of Lemma~\ref{l5.1} again gives a contradiction.
	Hence \(x_b\leq x_v\).
	Since \(x_b>0\) and
	\(x_v=\min_{w\in S}x_w\),
	it follows that
	\(
	x_a=\max_{z\in V(H)}x_z>0.
	\)
	
	We now prove that \(a\) is adjacent to every other
	vertex of \(H\).
	Suppose that \(az\notin E(H)\) for some
	\(z\in V(H)\setminus\{a\}\).
	Since \(S\) is a clique and \(ab\in E(H)\),
	we have \(z\notin S\cup\{b\}\).
	Since \(H\) is connected, choose \(y\in N_H(z)\)
	and replace the positive edge \(zy\) with \(za\).
	The negative triangle \(uabu\) remains.
Since \(x_a\geq x_y\) and \(x_a>0\), an argument
similar to that in Lemma~\ref{l5.1} gives a contradiction.
	Therefore \(N_H(a)=V(H)\setminus\{a\}\).

	We show that all coordinates on \(H\) are positive.
	For \(z\in V(H)\setminus\{a,b\}\), all edges incident
	with \(z\) are positive. Since \(az\in E(H)\), the
	eigenvalue equation at \(z\) gives
	\(
	\lambda x_z=\sum_{v\in N_P(z)}x_v
	\geq x_a>0.
	\)
	Together with \(x_a,x_b>0\), this proves that
	\(x_z>0\) for every \(z\in V(H)\).
	
	For \(z\in V(H)\setminus(S\cup\{b\})\), let
	\(\dot G_1'\) be obtained from \(\dot G^\ast\) by
	replacing the negative edge \(ub\) with the negative
	edge \(uz\).
	Clearly, \(\dot G_1'\) still has \(m\) edges and contains
	the negative triangle \(uazu\).
	It is also \(\mathcal K_r^-\)-free, since the degree
	of \(u\) remains at most \(r-2\).
	We have
	\(
	\Delta
	=\mathbf{x}^{\mathsf T}A(\dot G_1')\mathbf{x}
	-\mathbf{x}^{\mathsf T}A(\dot G^\ast)\mathbf{x}
	=2x_u(x_b-x_z).
	\)
	We consider the following two cases according to whether
	\(x_u=0\) or \(x_u>0\).
	
	If \(x_u=0\), the eigenvalue equation at \(u\) in
	\(\dot G^\ast\) gives
	\(
	x_b=\sum_{v\in S}x_v\geq x_a\geq x_b.
	\)
	Since all coordinates on \(H\) are positive, we obtain
	\(S=\{a\}\) and \(x_a=x_b\).
	Since \(x_u=0\), we have \(\Delta=0\), while
	\(
	\bigl((A(\dot G_1')-A(\dot G^\ast))\mathbf{x}\bigr)_u
	=x_b-x_z.
	\)
	By extremality and Lemma~\ref{lem:strict-increase},
	we must have \(x_z=x_b=x_a\), where \(z\in V(H)\setminus(S\cup\{b\})\).
		Thus all coordinates on \(H\) are equal.
		Since \(x_u=0\), the eigenvalue equations give
	\(d_H(v)=\lambda\) for every \(v\in V(H)\).
	Hence \(H\) is regular.
	Since \(a\) is adjacent to every other vertex of \(H\),
	it follows that \(H\) is complete.
Since \(S=\{a\}\), we have \(d_G(u)=2\), and hence
\(
p=m-2=e(H)=\binom{|V(H)|}{2}.
\)
Comparing this with \(p=\binom{k}{2}+s\), where
\(0\leq s\leq k-1\), we obtain \(|V(H)|=k\) and \(s=0\).
Therefore \begin{align}\label{1}
\dot G^\ast\cong\dot\Theta_{k,0}.
\end{align}

	If \(x_u>0\), we claim that 
	\(x_w>x_b\) for every \(w\in V(H)\setminus\{b\}.\)
For \(z\in V(H)\setminus(S\cup\{b\})\), if \(x_z\leq x_b\),
then
\(
\Delta=2x_u(x_b-x_z)\geq0.
\)
An argument similar to that in Lemma~\ref{l5.1}
gives a contradiction. Hence \(x_z>x_b\).
	For \(w\in S\), if \(x_w\leq x_b\), interchanging the
	signs of \(ub\) and \(uw\) gives a contradiction by
	a similar argument.
	Thus \(x_w>x_b\) for every \(w\in S\).
Therefore, \(b\) is the unique vertex
of \(H\) with minimum coordinate.
	We next show that for any distinct vertices
	\(v,w\in V(H)\) with \(x_v\geq x_w\), we have
	\(w\preceq v\), that is,
	\(
	N_H(w)\setminus\{v\}
	\subseteq N_H(v)\setminus\{w\}.
	\)
	Suppose to the contrary that
	\(z\in N_H(w)\setminus(N_H(v)\cup\{v\})\).
	Replace the positive edge \(wz\) with the positive edge
	\(vz\). This preserves the number of edges and the
	degree of \(u\).
	Since \(d_G(u)\leq r-2\) and \(ub\) remains the only
	negative edge, the resulting graph is
	\(\mathcal K_r^-\)-free.
	If \(wz\ne ab\), the negative triangle \(uabu\) remains.
	If \(wz=ab\), then
	\(w=a\), \(z=b\), and \(uavbu\) is a negative cycle
	in the resulting graph. Thus the graph is unbalanced.
	The Rayleigh change is
	\(
	\Delta=2x_z(x_v-x_w)\geq0.
	\)
	Since \(x_z>0\), an argument similar to that in
	Lemma~\ref{l5.1} gives a contradiction.
	Hence \(w\preceq v\).
	By Lemma~\ref{lem:threshold-preorder}, \(H\) is a threshold graph.
	
	We first consider the case that \(H\) is complete.
	Let \(H=K_c\). Since \(u\) has \(t\) positive neighbors
	and one negative neighbor in \(H\), we have
	\(m=\binom{c}{2}+t+1\) and \(p=m-2=\binom{c}{2}+t-1.\)
	The \(t+1\) neighbors of \(u\) are distinct, so
	\(1\leq t\leq c-1\), and hence \(0\leq t-1\leq c-2\).
	Comparing this with \(p=\binom{k}{2}+s\), where
	\(0\leq s\leq k-1\), we obtain \(c=k\) and \(t=s+1\).
	Therefore
	\begin{align}\label{2}
	\dot G^\ast\cong\dot B_{k,s+1}.
	\end{align}
	When \(s=0\), this graph is also \(\dot\Theta_{k,0}\).
	We next consider the case that \(H\) is not complete.
	
	\begin{lemma}\label{lem:noncomplete-case}
		If \(H\) is not complete, then \(s\geq1\) and
		\(\dot G^\ast\cong\dot\Theta_{k,s}\).
	\end{lemma}
	\begin{proof}[Proof of Lemma~\ref{lem:noncomplete-case}]
	
		Let \(R=\rho(R_p)\) and \(\eta=R-k+1\).
		By extremality and interlacing,
		\(
		\lambda\geq\lambda_1(\dot\Theta_{k,s})\geq R.
		\)
		If \(s=0\), then \(R=k-1\) and \(\eta=0\).
		If \(s>0\), then \(k-1<R<k\).
		Thus \(0\leq\eta<1\).
		By \eqref{a3.2}, we have
		\begin{equation}\label{eq:nc-R}
			\eta\bigl[R(R+1)-s\bigr]=s^2.
		\end{equation}

		Since \(H\) is a noncomplete threshold graph, choose
		a partition
		\(
		V(H)=C\mathbin{\dot\cup}I,
		\)
		where \(C\) is a maximum clique and \(I\ne\varnothing\)
		is an independent set.
		We may assume that \(b\in I\).
		Indeed, if \(b\in C\), choose \(z\in I\).
		Since \(x_z>x_b\), we have
		\(
		N_H(b)\setminus\{z\}
		\subseteq N_H(z)\setminus\{b\}.
		\)
		Thus \(z\) is adjacent to every vertex of
		\(C\setminus\{b\}\).
		The maximality of \(C\) gives \(zb\notin E(H)\).
		Moreover, \(b\) has no neighbor in \(I\), since such
		a neighbor would also be adjacent to \(z\), contrary
		to the independence of \(I\).
		Interchanging \(b\) and \(z\) therefore gives the
		required partition.
		Now \(S\subseteq N_H(b)\subseteq C\).
		For every \(z\in I\setminus\{b\}\), we have
		\(N_H(b)\subseteq N_H(z)\).
		Hence every vertex of \(S\) is adjacent to every other
		vertex of \(H\).

		Let \(T=\sum_{w\in V(H)}x_w\).
		For \(v\in S\), the eigenvalue equation gives
		\(
		(\lambda+1)x_v=T+x_u.
		\)
		Thus all coordinates on \(S\) are equal.
	By scaling \(\mathbf{x}\), we may assume that
	\(x_v=1\)
	for all \(v\in S.\)
		Since \(
		x_a=\max_{z\in V(H)}x_z
		\)
		and all coordinates on \(H\) are positive, we have
		\(0<x_w\leq1\) for every \(w\in V(H)\).
	Let \(|C|=L+1,\) \(I=\{z_1,\ldots,z_q=b\}\),
		\(d_i=d_H(z_i)\), and \(h=d_H(b)\).
	Since \(S\subseteq N_H(z_i)\) for every \(i\) and
	\(C\) is a maximum clique, after relabelling we have
	\[
	L\geq d_1\geq\cdots\geq d_q=h\geq t.
	\]
		Let \(D=\sum_{i=1}^q d_i\).
	Since \(p=m-2\), counting edges gives
		\begin{equation}\label{eq:nc-edges}
			D
			=m-\binom{L+1}{2}-t-1
			=p-\binom{L+1}{2}-(t-1).
		\end{equation}
		Since \(D\geq0\) and \(t\geq1\), we have
		\(
		p\geq\binom{L+1}{2}.
		\)
		Since \(p=\binom{k}{2}+s\), where \(0\leq s\leq k-1\),
		we have \(p<\binom{k+1}{2}\). Hence \(|C|=L+1\leq k\).

		By interlacing, \(\lambda\geq\rho(H)>L\).
		Deleting the negative edge \(ub\) strictly increases
		the Rayleigh quotient, since \(x_u,x_b>0\).
		Thus \(\lambda<\rho(P)\).
	Since \(e(P)=m-1=p+1\leq\binom{k+1}{2}\),
	Lemma~\ref{R} gives
	\(
	\rho(P)\leq\rho(K_{k+1})
	=k.
	\)
Hence
\begin{equation}\label{eq:lambda-range}
	L<\lambda<\rho(P)\leq k.
\end{equation}
	
		Since \(x_v\leq1\) for every \(v\in V(H)\), the eigenvalue
		equation at any vertex \(z_i\ne b\) gives
		\(
		\lambda x_{z_i}\leq d_i.
		\)
		Similarly, the equation at \(b\) gives
		\(
		\lambda x_b+x_u
		=\sum_{v\in N_H(b)}x_v\leq h.
		\)
		Since \(x_v=1\) for every \(v\in S\), the equation at \(u\)
		gives
		\(
		\lambda x_u=\sum_{v\in S}x_v-x_b=t-x_b.
		\)
		Combining the last two relations, 
		since \(\lambda>L\geq1\), we obtain
		\begin{equation}\label{eq:xb-bound}
			x_b\leq\frac{\lambda h-t}{\lambda^2-1}.
		\end{equation}
		
			Let
		\(\varepsilon=\lambda-L>0\), \(X_C=\sum_{v\in C}x_v,\)
		and \(X_I=\sum_{i=1}^q x_{z_i}.\)
	The eigenvalue equation at \(a\) gives
	\(
	\lambda+1=X_C+X_I+x_u.
	\)
	Summing the eigenvalue equations over \(C\), we obtain
	\(
	\varepsilon X_C
	=\sum_{i=1}^q d_ix_{z_i}+tx_u.
	\)
	Combining the above two equations, we obtain
	\(
	\varepsilon(\lambda+1)
	=\sum_{i=1}^q(d_i+\varepsilon)x_{z_i}
	+(t+\varepsilon)x_u.
	\)
		Multiplying the identity by \(\lambda\)
		and eliminating \(\lambda x_u\), we obtain
		\(
			\varepsilon\lambda(\lambda+1)
			\leq\sum_{i=1}^{q-1}d_i(d_i+\varepsilon)
			+t(t+\varepsilon)
			+
			\bigl[\lambda(h+\varepsilon)-(t+\varepsilon)\bigr]x_b.
		\)
	Since
	\(
	\lambda(h+\varepsilon)-(t+\varepsilon)
	\geq(\lambda-1)(t+\varepsilon)>0,
	\)
	we may substitute the upper bound for \(x_b\)
	given in \eqref{eq:xb-bound}.
		Simplifying gives
		\begin{equation}\label{eq:nc-master}
			(\lambda-L)\lambda(\lambda+1)
			\leq\sum_{i=1}^q d_i^2+(\lambda-L)D+\psi_t(\lambda,h),
		\end{equation}
		where, for \(x>L\),
		\[
		\psi_t(x,h)
		=\frac{(tx-h)
			\bigl[x^2+(t-1-L)x+L-h\bigr]}{x^2-1}.
		\]
		For \(t=1\), denote
		\begin{equation}\label{eq:nc-phi-def}
			\varphi(x,h)=\psi_1(x,h)
			=\frac{x-h}{x+1}
			\left(x-L+\frac{x-h}{x-1}\right).
		\end{equation}
		Since \(1\leq h\leq L<x\), we have
		\begin{equation}\label{eq:nc-phi}
			0<\varphi(x,h)<x-L+1.
		\end{equation}

		We next establish a numerical bound for the right-hand
		side of \eqref{eq:nc-master}.
		For \(t\leq h<L\), direct calculation gives
	\[
	\psi_t(x,h)-\psi_t(x,h+1)
	=\frac{x^2+(2t-1-L)x+L-2h-1}{x^2-1}
	\geq\frac{(x-L)(x+1)+(2t-2)x+1}{x^2-1}
	>0.
	\]
	Thus, \(\psi_t(x,h)>\psi_t(x,h+1)\).
	When \(t=1\), we have
	\begin{align}\label{a5.9}
	0<\varphi(x,h)-\varphi(x,h+1)<1.
	\end{align}
		
		Let \(\Delta=D+t-1\). By \eqref{eq:nc-edges}, we may write
		\(
		\Delta=p-\binom{L+1}{2}
		=QL+\nu,
		\)
		where \(Q\) and \(\nu\) are nonnegative integers, \(0\leq\nu<L\).
	Define \(\widehat{\nu}=\nu\) if \(\nu>0\), and \(\widehat{\nu}=L\) if \(\nu=0\).
	
	\begin{claim}\label{clm:nc-concentration}
		If \(L\geq2t-2\), then for every \(x>L\),
		\begin{equation}\label{eq:nc-concentration}
			\sum_{i=1}^q d_i^2+(x-L)D+\psi_t(x,h)
			\leq QL^2+\nu^2+(x-L)\Delta
			+\varphi(x,\widehat{\nu}).
		\end{equation}
	\end{claim}
	
	\begin{proof}[Proof of Claim~\ref{clm:nc-concentration}]
		We first reduce the numerical bound to \(t=1\),
		while keeping \(\Delta=D+t-1\) fixed.
		Suppose that \(t\geq2\).
		If \(h<L\), increase one smallest part from \(h\)
		to \(h+1\) and decrease \(t\) by one.
		This preserves \(\Delta\).
		The new minimum is at most \(h+1\), so the right-hand
		side of \eqref{eq:nc-master}, with \(\lambda\) replaced
		by \(x\), increases by at least
		\begin{equation}\label{eq:nc-shift}
			\begin{aligned}
				2h+1+x-L+\psi_{t-1}(x,h+1)-\psi_t(x,h)
				=\frac{2x(h-t+1)}{x-1}>2.
			\end{aligned}
		\end{equation}
		If all parts equal \(L\), add a part \(t-1\)
		and replace \(t\) by \(1\).
		This also preserves \(\Delta\), and the increase is
		\[
			(t-1)^2+(x-L)(t-1)
			+\varphi(x,t-1)-\psi_t(x,L)
			=\frac{x(L-2t+2)(x-L+1)}{x^2-1}\geq0.
		\]
		Thus, repeating these changes reduces the numerical
		bound to \(t=1\), when the sum of the parts is \(\Delta\).
	If two parts satisfy \(1\leq a_0\leq b_0<L\), replace
	them by \(a_0-1,b_0+1\).
	Let \(h'\) denote the minimum of the resulting positive
	parts, deleting the zero part when \(a_0=1\).
	The sum of the parts remains \(\Delta\), so the term
	\((x-L)\Delta\) is unchanged.
	If \(a_0\geq2\), then \(h'\leq h\), and hence
	\(\varphi(x,h')\geq\varphi(x,h)\).
	The increase in the sum of squares is
	\(
	(a_0-1)^2+(b_0+1)^2-a_0^2-b_0^2
	=2(b_0-a_0+1)>0.
	\)
	Thus the total change is
	\(
	2(b_0-a_0+1)+\varphi(x,h')-\varphi(x,h)>0.
	\)
	If \(a_0=1\), deleting the resulting zero part gives
	\(h=1\) and \(1\leq h'\leq b_0+1\). Therefore
	\(
	0\leq h'-h\leq b_0.
	\)
	The sum of squares increases by \(2b_0\).
By \eqref{a5.9}, we have
\(
\varphi(x,h)-\varphi(x,h')
\leq h'-h\leq b_0.
\)
	Consequently, the total change satisfies
	\(
	2b_0+\varphi(x,h')-\varphi(x,h)
	\geq b_0>0.
	\)
	Repeating these changes leaves at most one part smaller
	than \(L\).
	Since the sum remains \(\Delta=QL+\nu\), the final sequence
	consists of \(Q\) parts equal to \(L\), and one part
	\(\nu\) if \(\nu>0\).
	Its minimum is therefore \(\widehat{\nu}\), and
	\eqref{eq:nc-concentration} follows.
	\end{proof}
	
		\begin{claim}\label{clm:clique-size}
			\(|C|=k\).
		\end{claim}
		\begin{proof}[Proof of Claim~\ref{clm:clique-size}]
	Since \(|C|=L+1\leq k\), it suffices to prove that
	\(L=k-1\). Suppose to the contrary that \(L\leq k-2\).
		Let \(n_0=k-1-L\geq1\) and \(M_0=\binom{n_0+1}{2}\).
		
		\medskip
		\noindent\textbf{Case 1.}
		\textit{\(L\geq2t-2\).}
		
	Since \(\Delta=QL+\nu\), we have
	\(QL^2+\nu^2=L\Delta-\nu(L-\nu)\).
	Thus \eqref{eq:nc-master} and \eqref{eq:nc-concentration} give
		\begin{equation}\label{eq:nc-W}
			W(\lambda)\leq\varphi(\lambda,\widehat{\nu}),
		\end{equation}
		where
		\(
		W(x)=x\bigl[(x-L)(x+1)-\Delta\bigr]+\nu(L-\nu).
		\)
		For \(x\geq R\), we have
		\(x-L\geq1\)
		and
		\(x(x+1)\geq k(k-1)\geq p\geq\Delta.\)
		Hence
		\(
		W'(x)
		=(x-L)(2x+1)+x(x+1)-\Delta>1.
		\)
		On the other hand,
		differentiating \eqref{eq:nc-phi-def} with respect to \(x\),
		we obtain
		\(
		\frac{\partial\varphi}{\partial x}(x,h)
		=1-\frac{(h-1)^2}{2(x-1)^2}
		-\frac{2L(h+1)-h^2+1}{2(x+1)^2}<1,
		\)
		since \(1\leq h\leq L<x\).
		Thus \(W(x)-\varphi(x,\widehat{\nu})\) is strictly
		increasing for \(x\geq R\).

	Let
	\(\theta=\frac{s+\eta}{R}\leq1.\)
	Here 
	\(\theta\leq1\) follows from \(s\leq k-1\) and
	\(R=k-1+\eta\).
	Using \(k=L+n_0+1\), we obtain
	\(\Delta=\binom{k}{2}-\binom{L+1}{2}+s
	=n_0L+M_0+s\)
	and
	\(R=L+n_0+\eta.\)
By \eqref{eq:nc-R}, we obtain
\(
\eta(R+1)=s\theta.
\)
	Substituting the above identities into \(W(R)\)
	and simplifying, we obtain
	\begin{align}\label{eq:nc-WR}
	W(R)
	=(R-s)(M_0-n_0\theta)
	+s(M_0+s-L)+\nu(L-\nu).
	\end{align}

		\noindent\textbf{Subcase 1.1.} \textit{\(n_0\geq2\).}
		
		Let \(A_0=M_0-n_0\geq1\) and \(w_0=M_0+s-L\).
		Since \(\theta\leq1\),
		\(
		W(R)\geq(R-s)A_0+sw_0+\nu(L-\nu).
		\)
		Moreover, since \(0\leq\nu<L\), we have
		\(\nu(L-\nu)\geq0\).
		If \(w_0\geq A_0\), then
		\(
		W(R)\geq A_0R\geq R\geq n_0+\eta+1.
		\)
		If \(0\leq w_0<A_0\), then
		\(
		R-s=n_0+M_0+\eta-w_0>2n_0+\eta,
		\)
		so \(W(R)\geq(R-s)A_0>2n_0+\eta> n_0+\eta+1\).
		If \(w_0<0\), then \(\nu=M_0+s<L\) and
		\(
		sw_0+\nu(L-\nu)=M_0(L-\nu)>0.
		\)
		Since \(R-s>n_0+M_0+\eta\), 
		we have
		\(W(R)>(R-s)A_0> n_0+\eta+1\).
		By \eqref{eq:nc-phi},
		\(
		\varphi(R,\widehat{\nu})
		<R-L+1=n_0+\eta+1.
		\)
		Thus \(W(R)>\varphi(R,\widehat{\nu})\).
		
		\noindent\textbf{Subcase 1.2.} \textit{\(n_0=1\).}
		
	Since \(n_0=1\), we have \(L=k-2\), and hence
	\(s\leq k-1=L+1\). Let \(j=L+1-s\geq0\).
	Since \(R=\eta+k-1=L+1+\eta\) and \(j=L+1-s\), we have
	\(
	R-s=j+\eta.
	\)
		Then \(M_0=1\), \(1-\theta=\frac{j}{R}\), and
		\eqref{eq:nc-WR} gives
		\begin{align}\label{a11}
		W(R)=\frac{j(j+\eta)}{R}
		+s(2-j)+\nu(L-\nu).
		\end{align}
		
		If \(j\geq3\), then \(\nu=s+1<L\).
		Thus \(s(2-j)+\nu(L-\nu)=j-2\geq1\).
		Consequently,
		\(
		W(R)\geq\frac{j^2}{R}.
		\)
		Since \(R-\nu=j-1+\eta<j\) and \(R-L=1+\eta<2\),
		\eqref{eq:nc-phi-def} gives
		\(
		\varphi(R,\nu)=\frac{R-\nu}{R+1}
		\left(R-L+\frac{R-\nu}{R-1}\right)
		<\frac{3j}{R+1}<\frac{j^2}{R}\leq W(R).
		\)
		
		If \(j=2\), then \(s=L-1\) and
		\(\Delta=L+1+s=2L\), so \(\nu=0\) and
		\(\widehat{\nu}=L\).
		Since \(L\geq1\), we have
		\(R-L=1+\eta\leq R-1\).
		Thus \eqref{eq:nc-phi-def} gives
		\(
		\varphi(R,L)
		=\frac{1+\eta}{R+1}
		\left(1+\eta+\frac{1+\eta}{R-1}\right)
		\leq\frac{(1+\eta)(2+\eta)}{R+1}.
		\)
		Using \eqref{a11} and \(0\leq\eta<1\), we obtain
		\(
		W(R)=\frac{2(2+\eta)}{R}
		>\frac{2(2+\eta)}{R+1}
		>\frac{(1+\eta)(2+\eta)}{R+1}
		\geq\varphi(R,L).
		\)

	If \(j=0\), then \(s=L+1\).
	By \eqref{a11} and \(L\geq1\), we have
	\(
	W(R)=2s+\nu(L-\nu)
	\geq2(L+1)\geq4.
	\)
	Since \(R-L=1+\eta\) and \(0\leq\eta<1\),
	\eqref{eq:nc-phi} gives
	\(
	\varphi(R,\widehat{\nu})
	<R-L+1=2+\eta<3.
	\)
	Therefore \(W(R)>\varphi(R,\widehat{\nu})\).

If \(j=1\), then \(s=L\) and
	\(\Delta=2L+1\).
	If \(L\geq2\), then \(\nu=1\), and
	\eqref{a11} gives
	\(
	W(R)
	=2L-1+\frac{1+\eta}{R}
	>2L-1\geq3.
	\)
	On the other hand, \eqref{eq:nc-phi} gives
	\(
	\varphi(R,\widehat{\nu})
	<R-L+1=2+\eta<3.
	\)
	Hence \(W(R)>\varphi(R,\widehat{\nu}).\)
	If \(L=1\), then \(s=1\), \(\nu=0\), \(\widehat{\nu}=1\),
	and \(R=2+\eta\).
	Since \(s=1>0\), we have \(0<\eta<1\), so \(2<R<3\).
	Thus, \eqref{a11} gives
	\(
	W(R)=1+\frac{1+\eta}{R}
	=2-\frac1R
	>\frac{3}{2}.
	\)
	By \eqref{eq:nc-phi-def}, we have
	\(
	\varphi(R,1)=\frac{R(R-1)}{R+1}.
	\)
	Moreover,
	\(
	\frac32-\varphi(R,1)
	=\frac{(3-R)(2R+1)}{2(R+1)}>0.
	\)
	Consequently,
	\(
	W(R)>\frac32>\varphi(R,1)
	=\varphi(R,\widehat{\nu}).
	\)

		Thus \(W(R)>\varphi(R,\widehat{\nu})\) in every case.
		Since \(\lambda\geq R\) and \(W(x)-\varphi(x,\widehat{\nu})\) is strictly
		increasing for \(x\geq R\), we have
		\(
		W(\lambda)-\varphi(\lambda,\widehat{\nu})
		\geq W(R)-\varphi(R,\widehat{\nu})>0,
		\)
		contradicting \eqref{eq:nc-W}.
		Therefore \(L\leq k-2\) is impossible when \(L\geq2t-2\).

		\medskip
		\noindent\textbf{Case 2.}
		\textit{\(L\leq2t-3\).}

		Since \(d_i\leq L\),
		\(
		\sum_{i=1}^q d_i^2\leq LD.
		\)
		The monotonicity of \(\psi_t(x,h)\) in \(h\) gives
		\(
		\psi_t(x,h)\leq\psi_t(x,t)
		=t(x-L+t)\frac{x-1}{x+1}.
		\)
		Thus \eqref{eq:nc-master} yields
		\begin{equation}\label{eq:nc-E}
			E(\lambda)\leq0,
		\end{equation}
		where
		\begin{align}\label{a12}
		E(x)
		=x\bigl[(x-L)(x+1)-D\bigr]
		-t(x-L+t)\frac{x-1}{x+1}.
		\end{align}
	
Since \(L\geq h\geq t\), we have
\(
	E'(x)
	=3x^2+2(1-L)x-L-D-t
	+\frac{2t(L-t+1)}{(x+1)^2}
	>3x^2+2(1-L)x-L-D-t.
\)
By \eqref{eq:nc-edges} and \(p=\binom{k}{2}+s\), where
\(0\leq s\leq k-1\), we have
\(
	D+t
\leq p+1
	\leq
	\binom{k+1}{2}.
\)
For \(x\geq R\geq k-1\), using \(t\leq L\leq k-2\)
and \(
D+t
\leq
\binom{k+1}{2}
\), we obtain
\(
	3x^2+2(1-L)x-L-D-t
	\geq k(k+1)-1-(D+t)
	\geq\binom{k+1}{2}-1>0.
\)
Hence \(E'(x)>0\) for \(x\geq R\), so \(E(x)\) is strictly
increasing on \([R,\infty)\).

\noindent\textbf{Subcase 2.1.} \textit{\(s=0\).}

If \(s=0\), then \(\eta=0\), and hence
\(R=k-1=L+n_0\) and \(D=n_0L+M_0-t+1.\)
Substituting these identities into \eqref{a12}, we obtain
\(
	E(R)
	=R(M_0-1)+t(L-t)
	+\frac{2t(n_0+t)}{R+1}.
\)
Since \(n_0\geq1\), we have \(M_0\geq1\).
Together with \(L\geq t\), this shows that the first
two terms are nonnegative, while the last term is
strictly positive. Hence \(E(R)>0\).

		\noindent\textbf{Subcase 2.2.} \textit{\(s>0\).}
		
		Since \(t\leq L\leq2t-3\), we have \(t\geq3\).
		Moreover, \(L\leq k-2\) gives
		\(k\geq L+2\geq t+2\geq5\).
		Hence \(m=\binom{k}{2}+s+2\geq12\), and in particular
		\(m\geq7\).
		Since the range \(1\leq s\leq M\) has already been
		treated, we have \(s>M\).
		Thus either \(s=k-1\) or \(s\geq r-3\).
		Together with \(t\leq L\leq k-2\) and \(t\leq r-3\),
		this gives \(s\geq t\).
		
		Let
		\(d_0=s-t\geq0\)
		and
		\(f=\frac{R-1}{R+1}\in(0,1).\)
		Then \(D=n_0L+M_0+d_0+1\).
		Using \eqref{eq:nc-R} and \eqref{a12}, together with
		\(R-L=n_0+\eta\), we obtain
		\begin{align}\label{a13}
		E(R)
		=R(M_0-1-d_0)+s^2
		+\eta(n_0R+s-tf)-t(n_0+t)f.
		\end{align}
		
		Since \(s\geq t\), \(\eta\geq0\), and \(f<1\),
		\begin{equation}\label{eq:nc-ER}
			E(R)>
			R(M_0-1-d_0)+(2d_0-n_0)t+d_0^2.
		\end{equation}
		
		If \(n_0\geq2\) and \(M_0-1-d_0<0\), then
		\(
		R<L+n_0+1\leq2t+n_0-2.
		\)
		Substituting this upper bound into \eqref{eq:nc-ER}, we obtain
		\(
			E(R)
			>t(n_0^2-2)+(n_0-2)(M_0-1)
			+d_0(d_0-n_0+2)>0,
	\)
		where \(d_0>M_0-1\geq n_0\).
		
	If \(n_0\geq2\) and \(M_0-1-d_0\geq0\), then
	\(R=L+n_0+\eta\geq t+n_0\).
	Substituting this lower bound into \eqref{eq:nc-ER},
	we obtain
	\(
		E(R)
		>
		\frac{t(n_0-2)(n_0+1)}{2}
		+n_0(M_0-1)+(t-n_0)d_0+d_0^2.
	\)
	Since \(n_0\geq2\), the first term is nonnegative and
	\(M_0-1\geq n_0\).
	Moreover, since \(d_0\geq0\), we have
	\(
	d_0^2+(t-n_0)d_0
	=\left(d_0-\frac{n_0}{2}\right)^2
	+td_0-\frac{n_0^2}{4}
	\geq-\frac{n_0^2}{4}.
	\)
	Consequently,
	\(
	E(R)>n_0(M_0-1)-\frac{n_0^2}{4}
	\geq\frac{3n_0^2}{4}>0.
	\)
		
	If \(n_0=1\),
	then \(M_0=1\), and \(0\leq\eta<1\) gives
	\(
	R=L+1+\eta<L+2\leq2t-1.
	\)
Since \(s\geq t\) and \(f<1\), we have \(s-tf\geq0\).
Thus, by \eqref{eq:nc-R},
\(
\eta(R+s-tf)\geq\eta R
\geq\frac{s^2}{R+1}.
\)
	Consequently, by \eqref{a13}, we have
	\(
		E(R)
		\geq-d_0R+s^2-t(t+1)
		+\frac{s^2+2t(t+1)}{R+1}.
	\)
	Since \(d_0\geq0\), \(s=t+d_0\), and \(R<2t-1\),
	we further obtain
	\(
		E(R)
		>d_0^2+2d_0+1+\frac{t}{2}
		+\frac{d_0^2}{2t}>0.
	\)

		Thus \(E(R)>0\) in every case.
		Since \(E(x)\) is strictly
		increasing on \([R,\infty)\) and \(\lambda\geq R\),
		we have \(E(\lambda)\geq E(R)>0\),
		contradicting \eqref{eq:nc-E}.
		The assumption \(L\leq k-2\) is therefore impossible.
		Hence
		\(|C|=k\) and \(L=k-1.\)
	\end{proof}
We next prove that \(t=q=1\).
By Claim~\ref{clm:clique-size}, we have
\(L=|C|-1=k-1\).
Thus \eqref{eq:nc-edges} gives \(D+t-1=s\).
Moreover, \(s\leq k-1=L\).
	Since \(D\geq t\geq1\), we have
	\(s=D+t-1\geq2t-1\geq1\).
	Therefore, \eqref{eq:lambda-range} gives
	\(\lambda<k=L+1\).
	Thus, by \eqref{eq:nc-phi},
	\begin{align}\label{a14}
	\varphi(\lambda,h)<\lambda-L+1<2.
	\end{align}

	If \(t\geq2\), let \(F\) denote the right-hand side of
	\eqref{eq:nc-master}, that is,
	\(
	F=\sum_{i=1}^q d_i^2+(\lambda-L)D+\psi_t(\lambda,h).
	\)
	Since \(D=s-t+1\leq L-1\), we have \(h\leq D<L\).
	We may therefore apply the adjustment for \(t\geq2\)
	and \(h<L\) in the proof of Claim~\ref{clm:nc-concentration}.
	We increase one smallest part from \(h\) to \(h+1\)
	and decrease \(t\) by one.
	We repeat this adjustment until \(t=1\).
	Let \(F'\) denote the resulting value of \(F\).
	When \(t=1\), we have \(D=s=\sum_{i=1}^q d_i\), and hence
	\(
	\sum_{i=1}^q d_i^2
	\leq\left(\sum_{i=1}^q d_i\right)^2=s^2.
	\)
	Moreover, by \eqref{eq:nc-phi-def} and \eqref{a14}, we have
	\(\psi_t(\lambda,h)=\varphi(\lambda,h)<2\).
Therefore, \(F'<s^2+(\lambda-L)s+2\).
	By \eqref{eq:nc-shift}, we have
	\(
	F'-F>2.
	\)
	Thus
	\(F<s^2+(\lambda-L)s\).

		If \(t=1\) and \(q\geq2\), then \(D=s\geq2\) and
		\(
		\sum_{i=1}^q d_i^2
		\leq(s-1)^2+1\leq s^2-2.
		\)
		By \eqref{a14}, we obtain
		\(
			F
			\leq s^2-2+(\lambda-L)s+\varphi(\lambda,h)
			<s^2+(\lambda-L)s.
		\)
		
Let \(g(x)=(x-L)\bigl[x(x+1)-s\bigr]\).
Consequently, if \(t\geq2\) or \(q\geq2\), then
\eqref{eq:nc-master} implies \(g(\lambda)<s^2\).
Since \(g(x)\) is strictly increasing for \(x\geq R>L\)
and \(\lambda\geq R\), \eqref{eq:nc-R} gives
\(
g(\lambda)\geq g(R)
=\eta\bigl[R(R+1)-s\bigr]=s^2,
\)
a contradiction.
		Therefore
		\(t=q=1\) and \(d_1=s.\)

		Thus \(H\) consists of a positive \(K_k\) and the
		vertex \(b\), which has exactly \(s\) neighbors in
		the clique.
		Since \(S=\{a\}\subseteq N_H(b)\), the original signed
		graph is obtained by adding \(u\), the positive edge
		\(ua\), and the negative edge \(ub\).
		Hence \(\dot G^\ast\cong\dot\Theta_{k,s}\).
	\end{proof}
	
Combining \eqref{1}, \eqref{2}, and
Lemma~\ref{lem:noncomplete-case} with the reduction
at the beginning of this section, we conclude that
every extremal graph is switching isomorphic to
\(\dot\Theta_{k,s}\), or to \(\dot B_{k,s+1}\) when
\(0\leq s\leq M\).
If \(s=0\), then
\(\dot B_{k,1}\cong\dot\Theta_{k,0}\).
If \(s>M\), every extremal graph is switching isomorphic
to \(\dot\Theta_{k,s}\).
It therefore remains to compare the indices of the two
candidates when \(1\leq s\leq M\).

Let
\(\alpha=\lambda_1(\dot B_{k,s+1})\)
and
\(\gamma=\lambda_1(\dot\Theta_{k,s}).\)
	Recall that \(1\leq s\leq M\leq k-2\).
	By a direct calculation, the characteristic polynomials of
	\(\dot B_{k,s+1}\) and \(\dot\Theta_{k,s}\) are
	\((x+1)^{k-2}g_0(x)\) and \((x+1)^{k-3}f_0(x)\),
	respectively, where
	\(
		g_0(x)
		=x(x+1)(x-k+1)-(s+2)(x-k+1)-s^2
	\)
	and
	\(
	f_0(x)=x(x+1)g_0(x)-(k-3s)x-s(k-s-2).
	\)
	Since \((K_k,+)\subsetneq\dot B_{k,s+1}\)
	and the added vertex has signed degree \(s>0\),
	the Rayleigh principle gives \(\alpha>k-1\).
	Since \(\dot B_{k,s+1}\) is connected and unbalanced
	on \(k+1\) vertices, we also have \(\alpha<k\).
	Moreover, \(R_p\) is connected and \(K_k\subsetneq R_p\),
	so interlacing gives
	\(\gamma\geq\rho(R_p)>k-1\).
	Therefore, we have \(k-1<\alpha<k\)
	and
	\(\gamma>k-1.\)
Deleting the added vertex of \(\dot B_{k,s+1}\)
leaves a positive \(K_k\), so interlacing gives
\(\lambda_2(\dot B_{k,s+1})\leq\lambda_1(K_k)=k-1\).
Since \(s\leq k-2\), deleting \(a\) from
\(\dot\Theta_{k,s}\) leaves a signed graph with maximum
degree at most \(k-1\), so interlacing and the degree
bound give
\(
\lambda_2(\dot\Theta_{k,s})
\leq\lambda_1(\dot\Theta_{k,s}-a)
\leq\Delta(\dot\Theta_{k,s}-a)
\leq k-1.
\)
	Thus \(g_0\) and \(f_0\) have unique simple roots
	\(\alpha\) and \(\gamma\) in \((k-1,\infty)\), respectively.
	Since both polynomials are monic, for \(x>k-1\) we have
	\(\operatorname{sgn}g_0(x)=\operatorname{sgn}(x-\alpha)\)
	and
	\(\operatorname{sgn}f_0(x)=\operatorname{sgn}(x-\gamma).\)
	Since \(g_0(\alpha)=0\), we obtain
	\(
	f_0(\alpha)=-(k-3s)\alpha-s(k-s-2).
	\)
	It follows that
	\begin{equation}\label{eq:candidate-sign}
		\operatorname{sgn}(\gamma-\alpha)
		=\operatorname{sgn}\bigl((k-3s)\alpha+s(k-s-2)\bigr).
	\end{equation}
	
	If \(k\geq3s\), then both terms in
	\((k-3s)\alpha+s(k-s-2)\) are nonnegative.
	Their sum is zero only if
	\(k=3s\) and \(k=s+2\),
	that is, \((k,s)=(3,1)\).
	In this case,
	\(
	g_0(x)=(x-1)(x^2-5),
	\)
	so \(\alpha=\gamma=\sqrt5\), and both candidates are
	extremal.
	Otherwise, \eqref{eq:candidate-sign} gives
	\(\gamma>\alpha\), and the unique extremal graph is
	\(\dot\Theta_{k,s}\), up to switching isomorphism.
	
	Now suppose that \(k<3s\).
	Let
	\(d=3s-k>0\) and \(a_0=s(k-s-2)\geq0.\)
	Then \eqref{eq:candidate-sign} becomes
	\begin{align}\label{5.22}
	\operatorname{sgn}(\gamma-\alpha)
	=\operatorname{sgn}(a_0-d\alpha)
	=\operatorname{sgn}\left(\frac{a_0}{d}-\alpha\right).
	\end{align}
	
	If \(a_0\leq d(k-1)\), then
	\(
	\frac{a_0}{d}\leq k-1<\alpha,
	\)
	so \(\gamma<\alpha\).
	Hence the unique extremal graph is \(\dot B_{k,s+1}\),
	up to switching isomorphism.
	If \(a_0\geq dk\), then
	\(
	\frac{a_0}{d}\geq k>\alpha,
	\)
	so \(\gamma>\alpha\).
	Hence the unique extremal graph is \(\dot\Theta_{k,s}\),
	up to switching isomorphism.
	
	It remains to consider
	\(
	d(k-1)<a_0<dk,
	\)
	when
	\(
	k-1<\frac{a_0}{d}<k.
	\)
Using \(\operatorname{sgn}g_0(x)=\operatorname{sgn}(x-\alpha)\)
and \eqref{5.22}, we obtain
	\(
	\operatorname{sgn}(\gamma-\alpha)
	=\operatorname{sgn}\left(\frac{a_0}{d}-\alpha\right)
	=\operatorname{sgn}g_0\left(\frac{a_0}{d}\right).
	\)
	Multiplying by the positive integer \(d^3\), let
	\(
		N
		=d^3g_0\left(\frac{a_0}{d}\right)
		=a_0^3+(2-k)a_0^2d-(k+s+1)a_0d^2
		+\bigl[k(s+2)-s^2-s-2\bigr]d^3.
	\)
	Thus \(N\) is an integer and
	\(
	\operatorname{sgn}(\gamma-\alpha)=\operatorname{sgn}N.
	\)
	We claim that \(N\ne0\).
	Otherwise, \(g_0\left(\frac{a_0}{d}\right)=0\), so
	\(\frac{a_0}{d}\) is a rational root of the monic
	polynomial \(g_0(x)\) with integer coefficients.
	It must therefore be an integer, contradicting
	\(k-1<\frac{a_0}{d}<k\).
	Therefore, up to switching isomorphism, the unique
	extremal graph is \(\dot\Theta_{k,s}\) when \(N>0\),
	and \(\dot B_{k,s+1}\) when \(N<0\).
	This completes the proof of Theorem~\ref{thm:Kr-index}.
\end{proof}
\subsection{A partial characterization of the adjacency spectral radius extremal graphs}
	Throughout, let \(r\geq4\), \(m\geq3\), and 
	\(p=m-2=\binom{k}{2}+s,\) where \(0\leq s\leq k-1.\)
	Let \(\Lambda_{m,r}\) be the maximum index determined in
	Theorem~\ref{thm:Kr-index} and 
	\[
	\beta_{m,r}
	=\max\{\rho(H):H\text{ is an unsigned }K_r\text{-free graph},
	\ e(H)=m\}.
	\]
We first recall several spectral bounds for nonnegative
matrices and graphs.
	\begin{lemma}[\cite{HornJohnson2013}]
		\label{lem:row-sum-bound}
		Let \(B=(b_{ij})\) be an \(n\times n\) nonnegative
		matrix. Then
		\(
		\rho(B)\leq\max_{1\leq i\leq n}
		\sum_{j=1}^n b_{ij}.
		\)
	\end{lemma}
\begin{lemma}[\cite{Nikiforov2002}]\label{lem:Kr-unsigned-bound}
	Let \(r\geq3\) and 
	\(G\) be an unsigned \(K_r\)-free graph. Then
	\(
	\rho(G)^2\leq2e(G)\left(1-\frac{1}{r-1}\right).
	\)
\end{lemma}
\begin{lemma}[\cite{Nosal1970}]\label{lem:Nosal}
	Let \(G\) be an unsigned graph with \(m\) edges.
	If \(\rho(G)>\sqrt m\), then \(G\) contains a triangle.
\end{lemma}
The following lemma compares the index of a signed graph
with the spectral radius of its underlying graph.
In particular, the inequality is strict when the signed
graph is connected and unbalanced.
\begin{lemma}[\cite{Stanic2019Integral}]
	\label{lem:Kr-strict-underlying}
	Let \(\dot G=(G,\sigma)\) be a connected signed graph.
	Then \(\lambda_1(\dot G)\leq\rho(G)\), with equality
	if and only if \(\dot G\) is balanced.
\end{lemma}

		\begin{lemma}\label{lem:Kr-negative-end}
		Let \(\dot G\) be an unbalanced \(\mathcal K_r^-\)-free
		signed graph with \(m\) edges. Then
		\(
		-\lambda_n(\dot G)\leq\beta_{m,r},
		\)
		with equality if and only if \(\dot G\) is switching
		isomorphic to \((H,-)\), where \(H\) is an unsigned
		\(K_r\)-free graph with \(m\) edges and
		\(\rho(H)=\beta_{m,r}\).
	\end{lemma}
	
\begin{proof}
	Let \(\mu=-\lambda_n(\dot G)=\lambda_1(-\dot G)\).
	By Lemma~\ref{l1}, there exists a signed graph
	\(\widetilde{G}\) switching equivalent to
	\(-\dot G\) such that \(A(\widetilde{G})\) has
	a nonnegative unit eigenvector \(\mathbf{x}\)
	corresponding to \(\mu\).
	Let \(H\) be the spanning subgraph of \(\widetilde{G}\)
	consisting of the positive edges of \(\widetilde{G}\).
	We claim that \(H\) is \(K_r\)-free.
	Otherwise, an all-positive \(K_r\) in
	\(\widetilde{G}\) becomes an all-negative \(K_r\)
	in \(-\widetilde{G}\).
	Since \(-\widetilde{G}\) is switching equivalent
	to \(\dot G\), \(\dot G\) contains an unbalanced
	\(K_r\), a contradiction.
	By the Rayleigh principle,
	\(
		\mu
		\leq\mathbf{x}^{\mathsf T}A(H)\mathbf{x}
		\leq\rho(H).
	\)
	If \(e(H)<m\), choose a component \(F\) of \(H\) with
	\(\rho(F)=\rho(H)\), and attach \(m-e(H)\) pendant
	edges at a vertex of \(F\).
	The resulting graph \(H'\) is \(K_r\)-free and has
	\(m\) edges, while
	\(\rho(H')>\rho(F)=\rho(H)\).
	Consequently,
	\(
	\mu\leq\rho(H)<\rho(H')\leq\beta_{m,r}.
	\)
	If \(e(H)=m\), then
	\(\mu\leq\rho(H)\leq\beta_{m,r}\).
	Thus equality requires \(e(H)=m\), so
	\(\widetilde{G}\) is all-positive and \(\dot G\)
	is switching isomorphic to \((H,-)\), with
	\(\rho(H)=\beta_{m,r}\).

Conversely, let \(H\) be an unsigned \(K_r\)-free graph
with \(m\) edges and \(\rho(H)=\beta_{m,r}\).
We first show that \(H\) is connected.
Otherwise, choose a component \(F\) with
\(\rho(F)=\rho(H)\).
Delete the other components and attach \(m-e(F)\)
pendant edges at a vertex of \(F\).
The resulting graph is \(K_r\)-free and has \(m\) edges.
Since it is connected and properly contains \(F\), its
spectral radius is greater than \(\rho(F)\), a contradiction.
Thus \(H\) is connected.
We next show that \(\beta_{m,r}>\sqrt m\).
For \(a\geq1\), let \(\xi=\rho(K_{1,1,a})\).
The eigenvalue equations give \(\xi^2=\xi+2a\).
Since \(C_3\subseteq K_{1,1,a}\), we have
\(\xi\geq2\).
If \(m=2a+1\), then \(K_{1,1,a}\) is \(K_r\)-free
and has \(m\) edges, so
\(
\beta_{m,r}^2\geq\xi^2=\xi+2a>m.
\)
If \(m=2a+2\), attach a pendant edge to \(K_{1,1,a}\).
The resulting graph is \(K_r\)-free, has \(m\) edges,
and has spectral radius strictly greater than \(\xi\).
Since \(\xi^2=\xi+2a\geq m\), we again obtain
\(\beta_{m,r}>\sqrt m\).
Since \(\rho(H)=\beta_{m,r}>\sqrt m\),
Lemma~\ref{lem:Nosal} implies that \(H\) contains
a triangle.
Hence \((H,-)\) is unbalanced and
\(\mathcal K_r^-\)-free, and
\(
-\lambda_n((H,-))=\rho(H)=\beta_{m,r}.
\)
\end{proof}

	\begin{proof}[Proof of Theorem~\ref{thm:Kr-radius}]
	By Theorem~\ref{thm:Kr-index} and
	Lemma~\ref{lem:Kr-negative-end},
	\(
	\rho(\dot G)
	=\max\{\lambda_1(\dot G),-\lambda_n(\dot G)\}
	\leq\max\{\Lambda_{m,r},\beta_{m,r}\}.
	\)
The bound for \(\lambda_1\) is attained by an index
extremal graph, while the bound for \(-\lambda_n\)
is attained by \((H,-)\), where \(H\) is an unsigned
\(K_r\)-free graph with \(m\) edges and
\(\rho(H)=\beta_{m,r}\). Hence
	\begin{equation}\label{eq:Kr-radius-reduction}
		\max_{\dot G}\rho(\dot G)
		=\max\{\Lambda_{m,r},\beta_{m,r}\},
	\end{equation}
where the maximum is over all unbalanced
\(\mathcal K_r^-\)-free signed graphs with \(m\) edges.
Equality holds precisely for the index extremal graphs
when \(\Lambda_{m,r}\geq\beta_{m,r}\), and for
\((H,-)\) when \(\beta_{m,r}\geq\Lambda_{m,r}\),
where \(H\) is an unsigned \(K_r\)-free graph with
\(m\) edges and \(\rho(H)=\beta_{m,r}\).

			\medskip
		\noindent\textbf{Case 1.}
		\textit{\(3\leq m<\binom{r}{2}\).}
		
	Every unsigned graph with \(m\) edges is \(K_r\)-free.
	By Lemma~\ref{R},
	\(\beta_{m,r}=\rho(R_m)\), with \(R_m\) as the unique
	unsigned extremal graph.
	Let \(\dot G^\ast\) attain \(\Lambda_{m,r}\) and 
	\(G^\ast\) be its underlying graph.
	By Lemma~\ref{lem:one-component}, \(G^\ast\) is connected.
	Since \(\dot G^\ast\) is unbalanced and \(G^\ast\)
	is \(K_r\)-free, Lemma~\ref{lem:Kr-strict-underlying} gives
	\(
	\Lambda_{m,r}
	=\lambda_1(\dot G^\ast)
	<\rho(G^\ast)\leq\beta_{m,r}.
	\)
	Thus \eqref{eq:Kr-radius-reduction} gives
	\(\rho(\dot G)\leq\rho(R_m)\), with equality precisely
	for \((R_m,-)\).
		
		\medskip
		\noindent\textbf{Case 2.}
		\textit{\(m=\binom{r}{2}\).}
		
	An unsigned \(K_r\)-free graph with \(\binom{r}{2}\)
	edges is not isomorphic to \(K_r\).
	By Lemma~\ref{lem:second-extremal}, since \(L_r\)
	is \(K_r\)-free, we have
	\(\beta_{m,r}=\rho(L_r)\), with \(L_r\) as the unique
	unsigned extremal graph.
	Here
	\(p=\binom{r-1}{2}+r-3,\)
	\(k=r-1,\)
	and \(s=r-3>r-4.\)
	By Theorem~\ref{thm:Kr-index}, the unique index
	extremal graph is \(\dot\Theta_{r-1,r-3}\).
	Let \(F\) be its underlying graph.
	Since \(d_F(u)=2<r-1\) and \(F-u\cong R_p\) has
	\(r\) vertices and fewer than \(\binom{r}{2}\) edges,
	\(F\) is \(K_r\)-free.
	Since \(\dot\Theta_{r-1,r-3}\) is connected and unbalanced,
	Lemma~\ref{lem:Kr-strict-underlying} gives
	\(
	\Lambda_{m,r}
	=\lambda_1(\dot\Theta_{r-1,r-3})
	<\rho(F)\leq\beta_{m,r}.
	\)
	Therefore \(\rho(\dot G)\leq\rho(L_r)\), with equality
	precisely for \((L_r,-)\).

	\medskip
	\noindent\textbf{Case 3.}
	\textit{\(m\geq T_r:=\binom{2r-2}{2}+2\).}
	
	By Lemma~\ref{lem:Kr-unsigned-bound},
	\(
	\beta_{m,r}\leq
	U:=\sqrt{2m\left(1-\frac{1}{r-1}\right)}.
	\)
	Since \(\dot\Theta_{k,s}\) contains
	\((R_p,+)\) as an induced subgraph, interlacing gives
	\(\Lambda_{m,r}\geq\rho(R_p)\).
	If \(s>0\), let \(v\) be the vertex of \(R_p\) outside
	\(K_k\), and define a vector \(\mathbf{y}\) by
	\(y_w=1\) for every \(w\in V(K_k)\) and
	\(y_v=\frac{s}{k}\).
	Then
	\(
	\mathbf{y}^{\mathsf T}\mathbf{y}
	=k+\frac{s^2}{k^2}.
	\)
	Since \(d_{R_p}(v)=s\), we also have
	\(
	\mathbf{y}^{\mathsf T}A(R_p)\mathbf{y}
	=2\binom{k}{2}+2s\frac{s}{k}
	=k(k-1)+\frac{2s^2}{k}.
	\)
	By the Rayleigh principle,
	\(
		\rho(R_p)
		\geq
		\frac{\mathbf{y}^{\mathsf T}A(R_p)\mathbf{y}}
		{\mathbf{y}^{\mathsf T}\mathbf{y}}
		=k-1+\frac{s^2(k+1)}{k^3+s^2}
		\geq k-1+\frac{s^2}{k^2},
	\)
	where the last inequality follows from \(s^2\leq k^2\).
	The same bound holds for \(s=0\), since \(R_p=K_k\).
	Let \(z=\frac{s}{k}\). Thus
	\begin{equation}\label{eq:Kr-index-lower}
		\Lambda_{m,r}\geq\rho(R_p)\geq k-1+z^2.
	\end{equation}
	
Since \(m\geq T_r\) and \(s\leq k-1\), we have
\(
\binom{2r-2}{2}\leq p
=\binom{k}{2}+s<\binom{k+1}{2}.
\)
Hence \(k\geq2r-2\geq6\).
	By the definition of \(U\), it follows that
	\begin{align}\label{a6}
	U^2
	=2m\left(1-\frac{1}{r-1}\right)
	\leq2m\left(1-\frac{2}{k}\right).
	\end{align}
Since \(m=p+2=\binom{k}{2}+s+2\) and
\(z=\frac{s}{k}\), we have
\(2m=k(k-1)+2kz+4\).
By \eqref{eq:Kr-index-lower} and \eqref{a6}, we obtain
\[
\begin{aligned}
	\Lambda_{m,r}^2-U^2
	&\geq(k-1+z^2)^2
	-\bigl[k(k-1)+2kz+4\bigr]
	\left(1-\frac{2}{k}\right)\\
	&=k-5+\frac{8}{k}
	+2(k-1)z^2-2(k-2)z+z^4\\
	&\geq k-5+\frac{8}{k}
	-\frac{(k-2)^2}{2(k-1)}\\
	&=\frac{(k-4)^2}{2k}
	+\frac{k-2}{2(k-1)}>0.
\end{aligned}
\]
	Thus \(\Lambda_{m,r}>U\geq\beta_{m,r}\).
	By \eqref{eq:Kr-radius-reduction}, the maximum
	spectral radius is \(\Lambda_{m,r}\), with equality
	precisely for the index extremal graphs in
	Theorem~\ref{thm:Kr-index}.
		
			\medskip
	\noindent\textbf{Case 4.}
	\textit{\(\binom{r}{2}<m<T_r\).}
	
	Let \(H\) be a connected unsigned \(K_r\)-free graph
	with \(m\) edges and \(n_H\) vertices.
	For every \(v\in V(H)\), the row sum of \(A(H)^2\)
	at \(v\) is
	\(
		\sum_{w\in N_H(v)}d_H(w)
		\leq2m-n_H+1.
	\)
Applying Lemma~\ref{lem:row-sum-bound} to \(A(H)^2\),
we obtain
\(\rho(H)^2=\rho\bigl(A(H)^2\bigr)
\leq2m-n_H+1.\)
	Consequently, if \(\rho(H)\geq\Lambda_{m,r}\), then
	\begin{align}\label{a6.1}
	n_H\leq
	\left\lfloor2m+1-\Lambda_{m,r}^2\right\rfloor
	=N_{m,r}.
	\end{align}
	
Let \(\mathcal H_{m,r}\) be the family of all connected
unsigned \(K_r\)-free graphs with \(m\) edges and at most
\(N_{m,r}\) vertices.
If \(\beta_{m,r}\geq\Lambda_{m,r}\), every unsigned
extremal graph belongs to \(\mathcal H_{m,r}\).
If \(\beta_{m,r}<\Lambda_{m,r}\), then
\eqref{eq:Kr-radius-reduction} shows that the maximum
spectral radius is \(\Lambda_{m,r}\), attained precisely
by the index extremal graphs in
Theorem~\ref{thm:Kr-index}.
Therefore
\(
	\max\{\Lambda_{m,r},\beta_{m,r}\}
	=\max\left(
	\{\Lambda_{m,r}\}\cup
	\{\rho(H):H\in\mathcal H_{m,r}\}\right)
	=M_{m,r}.
\)
	By \eqref{eq:Kr-radius-reduction}, we obtain
	\(\rho(\dot G)\leq M_{m,r}\).
	Equality holds precisely for the index extremal
	graphs when \(\Lambda_{m,r}=M_{m,r}\), and for
	\((H,-)\), where \(H\in\mathcal H_{m,r}\) and
	\(\rho(H)=M_{m,r}\).
	
	Finally, we show that \(N_{m,r}\leq3r\).
	Since \(m<T_r=\binom{2r-2}{2}+2\), we have
	\(
	p=m-2<\binom{2r-2}{2}.
	\)
	Together with \(p\geq\binom{k}{2}\), this gives
	\(k\leq2r-3\).
	Let \(z=\frac{s}{k}\).
By \eqref{eq:Kr-index-lower}, \eqref{a6.1}, and
\(2m=k(k-1)+2kz+4\), we obtain
\[
\begin{aligned}
	N_{m,r}
	\leq
	\left\lfloor2m+1-(k-1+z^2)^2\right\rfloor
	=\left\lfloor
	k+4+2kz-2(k-1)z^2-z^4
	\right\rfloor
	\leq
	\left\lfloor
	k+4+\frac{k^2}{2(k-1)}
	\right\rfloor,
\end{aligned}
\]
where the last inequality follows from \(z^4\geq0\) and
\(
2kz-2(k-1)z^2
\leq\frac{k^2}{2(k-1)}.
\)
	Since \(k\leq2r-3\), we obtain
	\[
	\begin{aligned}
		k+4+\frac{k^2}{2(k-1)}
		=\frac{3k}{2}+\frac92+\frac{1}{2(k-1)}
		\leq3r+\frac{1}{2(k-1)}
		<3r+1.
	\end{aligned}
	\]
	Hence \(N_{m,r}\leq3r\).
	Thus \(M_{m,r}\) is obtained by comparing
	\(\Lambda_{m,r}\) with the spectral radii of a finite
	family of connected unsigned \(K_r\)-free graphs
	on at most \(3r\) vertices.
	This completes the proof of
	Theorem~\ref{thm:Kr-radius}.
	\end{proof}
	
	\section*{Declaration on the use of AI}
	
	The main proof ideas and overall strategies were developed
	by the authors. ChatGPT (OpenAI) was used to assist with
	literature searches, language editing, and LaTeX typesetting.
	The authors critically reviewed and revised all AI-assisted
	content, independently verified the mathematical statements
	and proofs, checked the references, and take full
	responsibility for the content of this paper.

\end{document}